\documentclass[review, preprint, authoryear]{elsarticle}
\usepackage[a4paper,left=2.54cm,right=2.54cm,top=2.54cm,bottom=2.54cm]{geometry}
\usepackage[bb=dsserif]{mathalpha}
\usepackage{bm}

\usepackage{amssymb}
\usepackage{amsthm}
\usepackage{amsmath}
\usepackage{blkarray}
\usepackage{amsfonts}
\usepackage{mathtools}
\usepackage{mathrsfs}
\usepackage{soul}
\setstcolor{blue}

\usepackage[normalem]{ulem}
\usepackage[dvipsnames]{xcolor}
\usepackage{tikz}
\usetikzlibrary{arrows.meta, backgrounds}
\usepackage{pgfplots}
\usepgfplotslibrary{statistics}
\usepackage{graphicx}
\usepackage[singlelinecheck=false, font=scriptsize, labelfont=scriptsize]{caption}
\usepackage{subcaption}
\usepackage{comment}
\usepackage{algorithm}
\usepackage{algpseudocode}
\algnewcommand\algorithmicforeach{\textbf{for each}}
\algdef{S}[FOR]{ForEach}[1]{\algorithmicforeach\ #1\ \algorithmicdo}
\algdef{S}[FOR]{ForAll}[1]{\algorithmicforeall\ #1\ \algorithmicdo}
\usepackage{enumitem}
\usepackage{booktabs}
\usepackage{natbib}
\usepackage{xurl}

\newcommand{\uh}[1]{\textcolor{black}{#1}}

\DeclareMathOperator*{\argmin}{arg\,min}

\journal{EURO Journal on Transportation and Logistics}

\begin{document}

\begin{frontmatter}



\title{Multi-period Stochastic Network Design for Combined Natural Gas and Hydrogen Distribution}


\author[a]{Umur Hasturk}
\author[b]{Albert H. Schrotenboer}
\author[a]{Kees Jan Roodbergen}
\author[a]{Evrim Ursavas}
\affiliation[a]{{Department of Operations, Faculty of Economics and Business, University of Groningen,},%Department and Organization
            addressline={Nettelbosje 2}, 
            city={Groningen},
            postcode={9747 AE}, 
            %state={},
            country={The Netherlands}}
\affiliation[b]{{Operations, Planning, Accounting and Control Group, School of Industrial Engineering, Eindhoven University of Technology,},%Department and Organization
            addressline={De Zaale}, 
            city={Eindhoven},
            postcode={5600 MB}, 
            %state={},
            country={The Netherlands}}

\begin{abstract}
Hydrogen is produced from water using renewable electricity. Unlike electricity, hydrogen can be stored in large quantities for long periods. This storage ability acts as a \textit{green battery}, allowing solar and wind energy to be generated and used at different times. As a result, green hydrogen plays a central role in facilitating a climate-neutral economy. However, the logistics for hydrogen are complex. As new pipelines are developed for hydrogen, there is a trend toward repurposing the natural gas network for hydrogen, due to its economic and environmental benefits. Yet, a rapid conversion could disrupt the balance of natural gas supply and demand. Furthermore, technical and economic developments surrounding the transition contribute additional complexity, which introduces uncertainty in future supply and demand levels for both commodities. To address these challenges, we introduce a multi-period stochastic network design problem for the transition of a natural gas pipeline network into a green hydrogen pipeline network. We develop a progressive-hedging-based metaheuristic to solve the problem. Results demonstrate our matheuristic is efficient, both in computation time and in solution quality. We show that factoring in uncertainty avoids premature expansion and ensures the development of an adequate pipeline network meeting long-term needs. In a case study in the Northern Netherlands for \textit{Hydrogen Energy Applications in Valley Environments for Northern Netherlands} initiative, we focus on two key scenarios: local production and importation, exploring their impacts on performance indicators. Our case insights exemplify the solid foundation for strategic decision-making in energy transitions through our approach.

\end{abstract}




\end{frontmatter}


\section{Introduction}

The European Union has set a goal to achieve carbon neutrality by 2050 \citep{eu2012}. This requires a greater reliance on renewable energy sources for electricity generation.  However, the supply from renewable energy sources, such as wind and solar, is intermittent and uncertain \citep{hodge2018combined}. Typically, natural gas power plants bridge the gap when demand surpasses renewable supply \citep{safari2019natural}. This approach, while practical, contradicts the carbon neutral directive and increases the dependency on countries rich in natural gas. An emerging and sustainable alternative is \textit{green hydrogen}, which is produced carbon-free through electrolysis powered by renewable sources \citep{megia2021hydrogen}. It allows the storage of excess energy as hydrogen, which can be converted back to power or directly used in transportation, residential heating, and industrial feedstock, highlighting its potential as a multifaceted energy solution \citep{staffell2019role}.

Today's energy landscape mainly centers around networks that consist of natural gas producers, customers, intermediate nodes, and storage units, linked by pipelines that facilitate the natural gas flow. In the coming decades, these networks are envisioned to gradually transition from natural gas to hydrogen-based infrastructure, a vision already taking shape with projects like \textit{Hydrogen Energy Applications in Valley Environments for Northern Netherlands} (HEAVENN) in the Netherlands, where Europe's first hydrogen valley is being built. This shift involves not only the construction of new hydrogen pipelines but also the retrofitting of existing natural gas pipelines for hydrogen distribution. Compared to building new pipelines, repurposing existing natural gas pipelines for hydrogen is both cost-effective and environmentally friendly \citep{gordon2023socio}. However, this process of repurposing is complex, as it requires a strategic approach to balance current natural gas supply and demand with the transition toward hydrogen to avoid potential disruptions. Rapid conversion could disrupt this balance, while delayed conversion could hinder the establishment of the hydrogen economy. Moreover, it is uncertain how fast hydrogen demand will emerge in society \citep{espegren2021role}. Consequently, it is crucial to strategically plan the construction of new pipelines and the gradual conversion of existing pipelines for the coming decades. 

In this study, we address a novel multi-period multi-commodity stochastic network design problem focused on the transition of gas-based pipelines into hydrogen-based pipelines over a given time horizon. At the beginning of the horizon, the pipeline network balances the existing supply and demand for all commodities (e.g., gas and hydrogen). As time progresses, we encounter a stochastic transition between these commodities, reflected in the changes in supply and demand at associated nodes in the network. Inventories are kept at each location, allowing for the storage of the commodities for the demand of the subsequent periods. We choose periods aligned with seasonality in the energy markets, which helps us to understand the interplay between inventory levels across different seasons. This approach utilizes hydrogen storage built up during summer's surplus to ensure adequate supply throughout the high-demand winter months. We employ a two-stage recourse model representation of this multi-period problem. The first stage involves making strategic decisions about the construction of new pipelines and retrofitting existing ones, while the second stage addresses uncertainty in demand and supply by adopting a scenario-based reformulation. Within each scenario, we determine the operational decisions related to commodity flow in pipelines and inventory management at network nodes. In cases where demand is not satisfied, we incur a penalty, mirroring real-world scenarios where trucks might be used temporarily to transport excess demand. 

To solve this problem, we design an efficient matheuristic based on progressive hedging. We augment this matheuristic with several acceleration techniques. \uh{Exact methods like Benders decomposition are also common for stochastic network design, but they face significant computational challenges with the multi-period complexity and scale of our problem, and are therefore not pursued in this study.} On well-calibrated benchmark instances, we show that our matheuristic outperforms a compact mixed-integer programming (MIP) formulation that is solved by Gurobi. Using the matheuristic, we study the structure of our optimal solutions. Among others, considering the stochastic nature of supply and demand leads to (i) an increased deployment of hydrogen infrastructure at the end of the planning horizon, reflecting a strategic investment that enhances the network's capacity and efficiency, (ii) a later adoption of hydrogen infrastructure, thereby reducing the risk of overcapacity and underutilization through better-informed planning, and (iii), a higher need for installation capacity in the peak adaptation periods as more pipelines need to be built and retrofitted in a shorter amount of time.

Our practice-based study makes several theoretical contributions to the theory of network design. Our problem relates to the literature on so-called stochastic fixed-charge capacitated multi-commodity network design problems \citep{hewitt2021stochastic}. Existing research focuses on balancing the uncertain supply and demand of multiple commodities by determining optimal arc construction ---pipelines in our context--- and flow distribution \citep{paraskevopoulos2016cycle, sarayloo2021learning}. In our paper, we offer several contributions to this literature. Firstly, we introduce a novel \textit{arc-sharing} setting within the stochastic fixed-charge capacitated multi-commodity network design problem, ensuring each pipeline can be used for any commodity, but is designated to one specific commodity and operates within its own distinct network in any particular period. %Secondly, while the existing literature primarily focuses on single-period analysis \citep{crainic2011progressive, yaghini2015cutting}, our problem addresses a multi-period framework, in line with the energy transition dynamics.
Secondly, we contribute to the area of incremental network design, where existing models address the multi-period recovery of partially disrupted networks \citep{engel2013incremental, kalinowski2015incremental}. While their multi-period strategic focus aligns with ours, their scope does not consider the multi-commodity aspect associated with the necessary transition between commodities over time. % one commodity network to another. 
Third, we combine elements of network design with storage and inventory decisions to reflect how seasonality in energy markets can be efficiently handled. Fourth, our enhancements to the progressive hedging algorithm extend its applicability beyond traditional network design to other two-stage recourse models across the related literature, offering a generalized framework that improves efficiency and scalability.
Finally, we present a case study in the Northern Netherlands associated with the HEAVENN initiative, which is a pioneering innovation in hydrogen economies \citep{newenergycoalition2020}. Our study focuses on two main scenarios with local production versus hydrogen import, and analyses how the transition evolves from natural gas to hydrogen. In summary, by evaluating how to effectively implement a transition from one commodity network to another over a time horizon, utilizing existing infrastructure, and constructing new pipelines when necessary, we provide valuable insights into the dynamics of the energy transition process, a perspective previously unexplored in the literature. 

The remainder of the paper is organized as follows. In Section \ref{sect:litrew}, we review the relevant literature. In Section \ref{sect:probDef}, we define the problem and the corresponding mathematical model. In Section \ref{sect:method}, we propose a progressive-hedging-based matheuristic to solve our problem. In Section \ref{sect:compExp}, we present and discuss the key insights and findings derived from our computational experiments. In Section \ref{sect:Case}, we tackle a case within the Northern Netherlands \citep{heavenncite}. In Section \ref{sect:concl}, we conclude our paper and present an outlook for future research. %In Appendix \ref{instGen}, we discuss about how we create the instances we text on Section \ref{sect:compExp}.

\section{Literature Review} \label{sect:litrew}
%\framebox{To be improved. More citations will be added.}

This section discusses several areas of the literature that relate to our work. In Section \ref{lite:hyd}, we consider studies with various aspects of hydrogen network design, encompassing supply chain planning, pipeline network configurations, and discussions on the optimal selection of pipeline diameters. Section \ref{lite:fixChar} discusses the literature on stochastic fixed-charge capacitated multi-commodity network design problems. We end this section by discussing how our work relates to existing studies on incremental network design.

\subsection{Hydrogen Network Design} \label{lite:hyd}

In the energy literature, several works focus on practical use cases within supply chain optimization. While effective in their specific applications, these models often provide limited theoretical insights from an Operations Research perspective. We highlight a selection of such works. \cite{almansoori2012design} present a multi-stage stochastic model for hydrogen supply chain network design. They encompass strategic decisions related to fueling stations, production plants, and storage facilities. Seeking to enhance the efficiency of such complex models, \cite{nunes2015design} employ sample average approximation and validate their approach with Great Britain’s liquid hydrogen infrastructure. We refer the interested reader to \cite{li2019hydrogen} for an optimization-oriented review of hydrogen supply chain network design.

Research on hydrogen pipeline network design has grown in recent years. For example, \cite{andre2013design} focus on constructing pipelines and their capacities for a case study in France. They incorporate nonlinear construction costs that vary depending on the determined diameters of the pipes. \cite{reuss2019modeling} assess the impact of linearized versus nonlinear diameter cost models on solution quality. A hydrogen production and transmission model is introduced by \cite{johnson2012spatially}. Applied to the Southwestern United States, they provide insights into supply location selection and the optimal dimensions of pipelines, showing its potential in effectively utilizing existing infrastructure. \cite{andre2014time} look at how hydrogen transport might change over time. Using a so-called backward heuristic, they predict pipeline growth in Northern France. Their work shows that trucks might be relevant for the short run, but pipelines become a better option when the hydrogen economy constitutes at least $10\%$ of the fueling market. This supports the idea in our paper of using a penalizing scheme, where trucks fill in when pipelines cannot meet demand for short-term imbalances.

To the best of our knowledge, prior studies have not explored the multi-commodity framework for converting existing pipelines from natural gas to hydrogen under a stochastic optimization setting. This includes the concept of arc-sharing, where not just two, but potentially multiple commodities are optimized simultaneously within their designated networks during the transition phase.
Moreover, unlike many existing models that focus on problem-specific, complex implementations, our model simplifies the problem structure by avoiding non-linear equations and excessive constraints, making it applicable for various pipeline systems. 
%Additionally, our advancements in the progressive hedging algorithm are applicable for other two-stage recourse models, creating a framework for the related literature. 
Furthermore, while previous studies often consider pipeline capacities as a continuous decision variable, our approach adopts a more realistic setting by selecting pipeline diameters from a limited set of standard options. Our study synthesizes these improvements to advance the current methodology in stochastic optimization for hydrogen network design.

%Furthermore, these studies consider pipeline capacities as a continuous decision variable, which does not reflect the industry practice of selecting from a limited set of standard pipeline diameters. Our approach incorporates a realistic approach by selecting pipeline diameters from these regional standard options.

%Moreover, unlike some prior models that restrict transportation modes to trucks and rails, our study specifically focuses to include pipelines with new constructions and repurposing option in the pipeline network. Additionally, the ones focusing on pipelines often model pipeline capacities as a continuous decision variable, our approach relates more closely to industry practices by selecting pipeline diameters from a limited set of standard options. %, providing a more practical framework for the deployment of pipeline networks. 
%By integrating these elements, our study provides a stochastic optimization framework that advances the current methodology applied in pipeline network design for hydrogen.}%, while also addressing the strategic planning under uncertainty. }

\subsection{Stochastic Fixed-Charge Capacitated Multi-commodity Network Design} \label{lite:fixChar}

The fixed-charge capacitated multi-commodity network design problem is an optimization problem that is faced in various domains including less-than-truckload (LTL) transportation, energy distribution, and telecommunications. In this design, the network distributes multiple commodities between their origin and destination via the selected capacitated arcs that incur a fixed cost (i.e., charge) and a variable cost per unit for commodity distribution. \cite{gray1971exact} introduces an exact method for a direct delivery case using decomposition and traditional branch-and-bound techniques. \cite{ghamlouche2003cycle} later address the general case, proposing a tabu search approach with a cycle-based neighborhood. However, real-world applications often lack complete upfront information ---such as demand and costs---, leading to the study of this problem under uncertainty. 

The literature concerning network design under uncertainty mostly considers two-stage recourse models to deal with demand uncertainty. The first stage determines arc capacity decisions, and the second stage decides the flow across chosen arcs, accounting for demand uncertainty via finitely many scenarios. Several solution techniques include a metaheuristic framework based on the progressive hedging algorithm \citep{crainic2011progressive}, an accelerated Bender's decomposition framework \citep{rahmaniani2018accelerating}, and partial Benders decomposition \citep{crainic2021partial} for the case where arc capacity is also uncertain. The literature also presents various extensions expanding upon the core concepts of this problem. \cite{hu2019multi} solve the problem with a multi-echelon setting, utilizing a multi-stage recourse model to address the problem. Considering the extreme scenarios, \cite{muller2021intermodal} propose fixing the first-stage decision variables with a high penalty incurred in the second stage. \cite{thapalia2012single} consider a variant of the problem where the arc capacity for a single commodity can be utilized in both directions of a directed network.

The chosen set of scenarios in the second stage significantly affects both the solution's quality and computational efficiency. In general, the more scenarios are considered, the closer the model is to reality. However, research has indicated that reducing the number of scenarios can accelerate solution times while maintaining the solution quality. In this context, \cite{crainic2014scenario} introduce an approach to grouping scenarios based on their similarities, thereby enhancing computational performance. Advancing the field, \cite{jiang2021soft} propose a soft clustering method, which, unlike traditional methods, allows scenarios to have probabilistic membership to multiple bundles.
 %This approach has demonstrated superior performance, achieving high solution quality in substantially reduced run-times.

Two-stage recourse models are also often selected to address the multi-period version of this design problem. Within the first stage, the decision is on which arcs should be selected throughout the entire planning horizon. Thereafter, the second stage decides on optimizing the commodity flow through these arcs for the planning horizon.  Adopting this paradigm, \cite{hewitt2019scheduled} and \cite{wang2019stochastic} partition an infinite time horizon into distinct, manageable periods, such as weeks, emphasizing periodic and cyclic scheduling. Taking a slightly different approach, \cite{jiang2021soft} and \cite{boland2017continuous} find solutions within a finite time horizon by constructing a time-space network. \cite{fragkos2021decomposition} presume that as the horizon progresses, arc construction becomes more economical, and present a solution to the deterministic version of the problem. 

Our contributions to this literature stream are significant. Firstly, in contrast to most of the existing literature where each commodity is often from a singular origin and a single destination, our study focuses on a small number of commodities dispersed across the network with multiple supply, demand, and storage nodes. This requires the model, contrary to the LTL setting, to uniquely assign each arc to a specific commodity, ensuring they operate independently within their own designated networks. Thus, one arc is only used for one commodity at a time in our arc-sharing context. Secondly, we incorporate inventory levels into the multi-period SND literature, enabling seasonal planning through strategic storage decisions for energy transition. Thirdly, our model addresses a transitional phase where the volume of at least one commodity diminishes over time. We leverage the existing network infrastructure of this diminishing commodity to benefit others for an efficient conversion of arcs from one commodity to another throughout the planning horizon. %To the best of our knowledge, none of these concepts have been previously considered within the stochastic network design optimization literature. 
To the best of our knowledge, the combination of these concepts have not been previously considered in this literature. 
Finally, we test this setting on a case study in the Northern Netherlands associated with the HEAVENN initiative, comparing scenarios of local production and hydrogen import to examine the transition from natural gas to hydrogen.

\subsection{Incremental Network Design} \label{lite:incND}

The incremental network design problem focuses on the restoration of networks impaired due to unforeseen events like natural disasters. Researchers aim to strategically rebuild or repair arcs over several periods. This concept is similar to our setting, where decisions on arc construction are made incrementally across a multi-period horizon. Constraints often include set budgets per period, repairs requiring more than one period, or a fixed number of arcs that may be built in each period.

\cite{kalinowski2015incremental} focus on this problem with an emphasis on maximizing the cumulative flow over the planning horizon by adding arcs in each period. Similar to our approach where additional arcs can be constructed every period, they operate with a single commodity from a single origin to a destination, aiming to optimize the flow of the commodity. On the other hand, \cite{baxter2014incremental} study a problem related to shortest paths in incremental network design. While they also consider the addition of arcs across periods, they restrict it to one arc per period, and focus on the minimization of the total length of flow prior to the complete repair of the network. 

Our research offers new perspectives on this existing literature. While we apply a multi-commodity framework to a multi-period network design problem ---diverging from traditional incremental network design--- we also introduce the arc-sharing concept. In our context, arc-sharing ensures that pipelines are designated for a single commodity, with the flexibility to switch commodities between periods, hence preventing cross-commodity interference and optimizing the flow within designated networks. Our approach, specifically tailored for the energy transition network design problem, therefore enhances the scope of incremental network design literature, especially in scenarios where diverse types of resources are considered as multiple commodities. Finally, our model incorporates inventory management via storage decisions, enabling strategic resource planning not only for the current period but also for future periods, addressing long-term resource allocation.
 
\section{Problem Definition} \label{sect:probDef}
We first provide an overview of our system and detail the dynamics within. Afterward, we introduce the MIP that models our energy transition network design problem.

\subsection{System Description}
We consider an energy transition network design problem on a finite planning horizon with a directed network $\mathscr{G} = (\mathscr{N},\mathscr{A})$, where $\mathscr{N}$ is the node set and $\mathscr{A}$ the arc set. We consider a setting with periods $t \in \mathscr{T} \coloneqq \{ 0, 1,\ldots,T - 1\}$. \uh{At each time period $t$, decisions regarding infrastructure (construction or conversion) can be implemented and network flows are to be determined. For the purpose of this study, we define each period as a six-month interval (e.g., a winter or summer season), a discretization chosen to capture seasonal variations in energy supply and demand. However, alternative time discretizations can be adopted without altering the model.} 

At period $t=0$, we assume the network is constructed so that it enables matching demand and supply in the existing energy market. For this, we consider multiple commodities $k \in \mathscr{K} \coloneqq \{ 1,2,\ldots, K\}$ that model the different energy types\uh{, such as natural gas or hydrogen}. Each node $i \in \mathscr{N} \coloneqq \{ 1,2,\ldots,N \}$ has a net supply quantity for each commodity, which is negative for demand nodes and zero for intermediate nodes. \uh{The demand may be satisfied by any supplier within the same commodity network, with no fixed origin-destination pair as seen in classical multi-commodity network design. This assumption reflects a national-level energy system without contract-based routing. If needed, the model can accommodate more detailed commodity definitions, such as separate purity grades of hydrogen for distinct use in industrial or transportation applications, each represented as a distinct commodity.} For some commodities, there may be no market initially, which means that their net supply quantities are zero for each node. This also means that an intermediate node for a commodity could become a supply or a demand node in later periods of the time horizon, for example when an electrolysis facility for hydrogen is built at that node. We assume each of these nodes is included in the initial set $\mathscr{N}$, and thus the node set does not change over periods of the horizon. One node may be shared by multiple commodities, for example, a supplier node for one commodity may be a demand node for another commodity. 

The arc set $\mathscr{A} \coloneqq \{ 1,2,\ldots,A\}$ consists of all initial arcs and all candidate arcs to be constructed. There may be multiple arcs between a pair of nodes, for example, each representing a pipeline with a different standard capacity. $\mathscr{A}_k \subseteq \mathscr{A}$ is the set of arcs that are in use at $t=0$ for commodity $k \in \mathscr{K}$, where $\mathscr{A}_i \cap \mathscr{A}_j = \varnothing$. $\mathscr{A}_0$ represents the candidate arcs that are not constructed at $t=0$ and thus might be constructed in future periods. The cardinality of the set of arcs $\mathscr{A}$, thus, does not change over periods of the horizon. Once an arc is constructed, it is available in all the following periods. An arc can only be assigned to a single commodity $k$ at any given time, but this assignment can switch to a different commodity during the planning horizon. We enable the switch between commodities via decision variables, detailed below, and the set $\mathscr{A}_k$ does not change over time. This conversion process is less costly compared to constructing a new arc. We impose that each arc can be converted at most once during the planning horizon. %Furthermore, we assume construction and conversion occur in between periods, thus once constructed, an arc is available on all the following periods.

As is also common in multi-period stochastic network design literature, we assume no control over location decisions \citep{boland2017continuous, fragkos2021decomposition, jiang2021soft}. Thus, we model net supply at each node for each period as a stochastic parameter, capturing both periodic variations, such as seasonal effects, and the uncertainty inherent in the emerging hydrogen economy. Specifically, it is uncertain when and where supplier nodes will emerge as the hydrogen economy develops. This uncertainty naturally extends to storage capacities at nodes, which depend on their roles within the supply chain. Supplier nodes are typically associated with larger, long-term storage units, whereas demand nodes, such as residential areas or refueling stations, require smaller, short-term storage capacities. We reflect this correlation in our instance generator, and assign suppliers and their inventory capacities accordingly. One could argue that incorporating supply and storage capacities as decision variables could enhance the study. However, this would require a broader scope integrating network expansion planning. In addition, it can be expected, as we currently observe in the oil, fuel, and gas industry, that storage solutions will be exploited by other parties.  Such an extension would thus require integrating multiple objectives, perspectives, and incentives into our model. This would shift the focus away from our main contribution: understanding the transition of natural gas pipelines to hydrogen pipelines. 
%and thus treat both supply decisions and storage capacities as stochastic parameters.

\uh{We model this problem as a two-stage mixed-integer program with continuous recourse. This is consistent with the complexity observed in stochastic network design problems, where two-stage approximations are commonly used to replace multi-stage formulations \citep{wang2019stochastic, jiang2021soft}. Similar simplifications are found in scenario tree generation \citep{heitsch2009scenario} and rolling-horizon methods \citep{cavagnini2022rolling}, which often emulate multi-stage dynamics through tractable two-stage constructs. More importantly, high-pressure transmission pipelines, whether newly constructed or repurposed, typically require 24–60 months from final investment decision to commissioning \citep{ENTSOG2024}. Such long lead times limit the practicality of multi-stage recourses and thus favor a two-stage representation in which infrastructure commitments are fixed prior to the planning horizon \citep{beste2019adaptive}. We assume that the model is implemented with sufficient advance notice to allow all preparatory activities to be completed, ensuring that construction can proceed as scheduled within the planning horizon. If needed, the two-stage model can always be run in a rolling-horizon fashion, but given the long lead time between planning and execution, we limit our scope to the canonical two-stage mixed-integer program with continuous recourse.}
%More importantly, a two-stage approach is often closer to reality due to the extended throughput time from conception to realization, typically spanning several years, which limits the feasibility of making frequent changes. 
%In line with these, our two-stage model implements a finite set of scenarios at the second stage. 

In the first stage, we decide on which arcs to construct for the entire planning horizon, and simultaneously assign a specific commodity to each newly constructed arc, as well as decide on converting existing arcs from one commodity to another. In the second stage, we model the uncertainty using a finite set of scenarios. We optimize for each scenario the distribution of each commodity to ensure supply and demand are balanced. 
We further incorporate in our model the flexibility to handle unmet demand by penalizing it, simulating alternative transport options such as trucks.
To address practical questions regarding seasonal storage, our model adopts a two-season framework per year ---winter and summer--- as commonly practiced in the HEAVENN region, where, for example, excess supply is allowed to be stored and carried over to the next seasons. The details of this structure, as well as how we generate our instance parameters in our study are given in  \ref{instGen}. The goal is to find a network design that minimizes the construction cost and the expected distribution cost. % allowing storing and utilizing the products in subsequent periods.we account for uncertainty by representing it through finitely many scenarios. For each scenario, we determine the inventory flow on constructed arc network, and the amount of inventory to be stored at locations for future use.

%The primary objective of this problem is to make decisions for each time period on which arcs to construct and which arcs to convert from one commodity to another.  %The aim of the problem is to decide for each periods which arcs to construct, which arcs to convert from one commodity to another, and how to distribute each commodity with the available arc network, and how much to store commodities on nodes to balance supply and demand for all commodities. 

 %for some extreme cases. 
 %, and when to pay penalty for unsatisfied amounts of demand.

\subsection{Mixed Integer Programming Formulation}
Let $ y_{at}^k $ be a binary decision variable, equaling $1$ if pipeline arc $a \in \mathscr{A}$ is constructed at time period $t$ for commodity $k$, which has a fixed construction cost $ f_{at}^k $. Let $ b_{at}^k $ be a binary decision variable, equaling $1$ if pipeline arc $a$ is converted for commodity $k$ at the beginning of time period $t$, which has a fixed conversion cost $ g_{at}^k$. Each arc $a$ has a fixed flow capacity of $U_a$. Let $ z_{at}^k $ be a binary decision variable, having a value of $1$ if the pipeline arc $a$ is assigned to be used at time period $t$ for commodity $k$. 

We model uncertainty with finitely many scenarios $\omega \in \Omega \coloneqq \{ 1,2,\ldots,W\}$. Each scenario $\omega$ has a probability of occurrence of $ p_\omega$. Under scenario $\omega$, $ s_{it}^{k\omega}$ represents the net supply of node $i$ for commodity $k$ in period $t$. Let $ x_{at}^{k\omega} $ be the flow amount of commodity $k$ on arc $a$ on period $t$ under scenario $\omega$, which costs $c_{at}^{k}$ per unit, which we set linearly proportional to the arc length. Let $I_{it}^{k\omega}$ be the inventory level of node $i$ for commodity $k$ at the beginning of period $t \in \mathscr{T} \cup \{T\}$ under scenario $\omega$, which cannot exceed the inventory capacity of $C_{it}^{k\omega}$. \uh{Consistent with $s_{it}^{k\omega}$, we treat $C_{it}^{k\omega}$ as a scenario-dependent parameter, reflecting exogenous uncertainty in future infrastructure deployment. Capacity levels depend on the evolving role of each location in the hydrogen economy; for example, future suppliers are more likely to be assigned larger storage units.}
Holding costs, which typically include depreciation, risk of damage, and opportunity costs \citep{FALLAH2011199}, are less accounted for in hydrogen applications \citep{heavenncite}. Thus, we assume no inventory holding costs in this study.  However, the model can be adjusted to include these costs if necessary, without affecting the methodology presented in this paper.
Let $e_{it}^{k\omega}$ be the amount of unsatisfied net demand in node $i$ for commodity $k$ in period $t$ under scenario $\omega$, which is penalized with a cost of $h$ per unit. This penalty may, for example, represent the variable cost of emergency transport by truck. 
Let the sets $\delta^+(i)$ and $\delta^-(i)$ represent the arc sets that enter and leave node $i$, respectively.

The two-stage stochastic mixed integer program with continuous recourse is formulated as:
\allowdisplaybreaks
\begin{align}
    \text{min  } & \sum_{a \in \mathscr{A}} \sum_{k \in \mathscr{K}} \sum_{t \in \mathscr{T}} \left( f_{at}^k y_{at}^k  + g_{at}^k b_{at}^k \right) + \sum_{\omega \in \Omega} p_\omega \sum_{k \in \mathscr{K}} \sum_{t \in \mathscr{T}} \left(  \sum_{a \in \mathscr{A}}   c_{at}^k x_{at}^{k\omega}  + \sum_{i \in \mathscr{N}}  h e_{it}^{k\omega} \right) \span \label{firstEq}\\
    \text{s.t  } & I_{i,t+1}^{k\omega} \leq I_{it}^{k\omega} + s_{it}^{k\omega} + e_{it}^{k\omega}- \sum_{a \in \delta^-(i)} x_{at}^{k\omega} + \sum_{a \in \delta^+(i)} x_{at}^{k\omega} \hspace{0.1cm} & \forall i \in \mathscr{N}, k \in \mathscr{K}, \omega \in \Omega, t \in \mathscr{T}  \label{balance} ,\\
    & I_{it}^{k\omega} \leq C_{it}^{k\omega} & \forall i \in \mathscr{N}, k \in \mathscr{K}, \omega \in \Omega, t \in \mathscr{T} \cup \{T\} \label{invCapacity},\\
    & x_{at}^{k\omega} \leq U_a  z_{at}^{k} & \forall a \in \mathscr{A}, k \in \mathscr{K}, \omega \in \Omega, t \in \mathscr{T} \label{arcCapacity},\\
    & \sum_{k \in \mathscr{K}} \sum_{t \in \mathscr{T}} y_{at}^k \leq 1 & \forall a \in \mathscr{A} \label{arcConstruction} ,\\
    & y_{at}^{k} \leq z_{at}^{k} & \forall a \in \mathscr{A}, k \in \mathscr{K}, t \in \mathscr{T} \label{arcAssignFirst},\\
    & \sum_{k \in \mathscr{K}} \sum_{t' = 0}^t y_{at'}^{k} = \sum_{k \in \mathscr{K}} z_{at}^{k}  & \forall a \in \mathscr{A}, t \in \mathscr{T} \label{arcAssignLater},\\
    & z_{at}^{k} + \sum_{k' \in \mathscr{K} \setminus \{k\}} z_{a,t-1}^{k'} \leq 1 + b_{at}^k & \forall a \in \mathscr{A}, k \in \mathscr{K}, t \in \mathscr{T} - \{0\} \label{conversion},\\ 
    & \sum_{k \in \mathscr{K}} \sum_{t \in \mathscr{T}} b_{at}^k \leq 1 & \forall a \in \mathscr{A} \label{conversionLimit} ,\\
    & y_{a0}^k = 1 & \forall k \in \mathscr{K}, a \in \mathscr{A}_k \label{initialNetwork} ,\\ 
    & I_{i0}^{k\omega} = o_{i}^{k} &  \forall i \in \mathscr{N}, k \in \mathscr{K}, \omega \in \Omega \label{initInv}, \\
    & x_{at}^{k\omega} \geq 0 & \forall a \in \mathscr{A}, k \in \mathscr{K}, \omega \in \Omega, t \in \mathscr{T},\\
    & I_{it}^{k\omega} \geq 0 & \forall i \in \mathscr{N}, k \in \mathscr{K}, \omega \in \Omega, t \in \mathscr{T} \cup \{T\}, \\
    & e_{it}^{k\omega} \geq 0 & \forall i \in \mathscr{N}, k \in \mathscr{K}, \omega \in \Omega, t \in \mathscr{T},\\
    & y_{at}^k \in \{0,1\} & \forall a \in \mathscr{A}, k \in \mathscr{K}, t \in \mathscr{T},\\
    & z_{at}^{k}\in \{0,1\} & \forall a \in \mathscr{A}, k \in \mathscr{K}, t \in \mathscr{T},\\
    & b_{at}^{k}\in \{0,1\} & \forall a \in \mathscr{A}, k \in \mathscr{K}, t \in \mathscr{T} .\label{lastEq}
\end{align}

The objective function minimizes the total cost of arc construction, arc conversion, the expected cost of arc flow, and the expected penalty cost. Constraint \eqref{balance} ensures the inventory levels are balanced between periods according to the net supply, as well as the inflow and outflow of the commodity. 
Constraint \eqref{invCapacity} ensures that the inventory respects the capacity. Constraint \eqref{arcCapacity} ensures that the flow capacity of an arc is not exceeded, and is only available for the assigned commodity. Constraint \eqref{arcConstruction} ensures each arc can be built once in the planning horizon. 
Constraint \eqref{arcAssignFirst} ensures that at the period an arc is constructed, it is assigned to its corresponding commodity. In the following periods, the arc may be converted for the use of another commodity. For those periods, Constraint \eqref{arcAssignLater} satisfies that each constructed arc is assigned to one of the commodities.  
We link the variable $z$ with the conversion variable $b$ by Constraint \eqref{conversion}, ensuring the conversion cost is reflected in the objective function value. Each arc is restricted to be converted at most once during the planning horizon by Constraint \eqref{conversionLimit}. The initial networks for each commodity are created by Constraint \eqref{initialNetwork}. Finally, with Constraint \eqref{initInv} we make sure that the initial inventories are assigned to the start inventory values of $o_{i}^{k}$. The remaining constraints model the decision variables' domains. \uh{An overview of the notation is given in Table \ref{tab:NOTATION}.}

\begin{table}[t]
\centering
 
%\captionsetup{ labelfont={color=red}, textfont={color=red} }
\caption{Overview of notation} \label{tab:NOTATION}
\begin{tabular}[h]{p{0.08\linewidth}  p{0.88\linewidth}} 
\toprule
\multicolumn{2}{l}{Parameters:} \\
$f_{at}^k$ & fixed construction cost of pipeline $a$ at period $t$ for commodity $k$ \\
$g_{at}^k$ & fixed conversion cost of pipeline $a$ at period $t$ for commodity $k$ \\
$p_{\omega}$ & probability of occurrence of scenario $\omega$ \\
$c_{at}^k$ & unit flow cost of pipeline $a$ at period $t$ for commodity $k$ \\
$h$ & unit penalty cost of unsatisfied demand \\
$s_{it}^{k\omega}$ & net supply of node $i$ at period $t$ for commodity $k$ under scenario $\omega$ \\
$C_{it}^{k\omega}$ & inventory capacity of node $i$ at period $t$ for commodity $k$ under scenario $\omega$ \\
$U_a$ & fixed flow capacity of pipeline $a$ \\
\multicolumn{2}{l}{Decision variables:} \\
$y_{at}^k$ & binary decision variable equaling $1$ if pipeline $a$ is constructed at period $t$ for commodity $k$ \\
$b_{at}^k$ & binary decision variable equaling $1$ if pipeline $a$ is converted at period $t$ for commodity $k$ \\
$z_{at}^k$ & binary decision variable equaling $1$ if pipeline $a$ is assigned at period $t$ for commodity $k$ \\
$x_{at}^{k\omega}$ & flow amount of pipeline $a$ at period $t$ for commodity $k$ under scenario $\omega$ \\
$e_{it}^{k\omega}$ & unsatisfied demand of node $i$ at period $t$ for commodity $k$ under scenario $\omega$ \\
$I_{it}^{k\omega}$ & inventory level of node $i$ at the beginning of period $t$ for commodity $k$ under scenario $\omega$ \\
\bottomrule     
\end{tabular}
\end{table}

\section{Progressive Hedging Based Matheuristic} \label{sect:method}

The progressive hedging algorithm (PHA) is a method that uses an iterative process to solve two-stage recourse optimization problems, such as stochastic network design \citep{crainic2011progressive, crainic2014scenario, jiang2021soft, sarayloo2023integrated}. By implementing the Augmented Lagrangian method \citep{jia2024scenario}, the problem can be decomposed among each scenario. Then, the PHA progresses iteratively by modifying the objective function of the single scenario deterministic optimization problems. The modification is done in such a way that the solutions to the deterministic optimization problems become closer to each other as the algorithm progresses through these iterations. The algorithm continues until the solutions of each single scenario problem share the same first-stage solution variables. 

To apply the PHA to our model, we first create a reformulation of model \eqref{firstEq} - \eqref{lastEq}. This reformulation duplicates each first-stage variable for each scenario, and then adds non-anticipativity constraints to make sure these variables are equal to each other. Explicitly, we create a copy for each first-stage variable for each scenario. We denote these variables by $y_{at}^{k\omega}$, $z_{at}^{k\omega}$, and $b_{at}^{k\omega}$. % for all $a \in \mathscr{A}, k \in \mathscr{K}, \omega \in \Omega, t \in \mathscr{T}$. 
Furthermore, we add the non-anticipativity Constraints \eqref{nac1} - \eqref{nac3}. This results in the following reformulation:% which would yield the same optimal value of the model \eqref{firstEq} - \eqref{lastEq}.
\allowdisplaybreaks
\begin{align}
    \text{min  } & \sum_{\omega \in \Omega} p_\omega \Bigg[ \sum_{a \in \mathscr{A}} \sum_{k \in \mathscr{K}} \sum_{t \in \mathscr{T}} \left( f_{at}^k y_{at}^{k\omega}  + g_{at}^k b_{at}^{k\omega} \right) +  \sum_{k \in \mathscr{K}} \sum_{t \in \mathscr{T}} \left(  \sum_{a \in \mathscr{A}}   c_{at}^k x_{at}^{k\omega}  + \sum_{i \in \mathscr{N}}  h e_{it}^{k\omega}\right) \Bigg] \span  \label{phStart}\\
    \text{s.t  } & I_{i,t+1}^{k\omega} \leq I_{it}^{k\omega} + s_{it}^{k\omega} + e_{it}^{k\omega} - \sum_{a \in \delta^-(i)} x_{at}^{k\omega} + \sum_{a \in \delta^+(i)} x_{at}^{k\omega} \hspace{0.1cm} & \forall i \in \mathscr{N}, k \in \mathscr{K}, \omega \in \Omega, t \in \mathscr{T}  \label{balancePH}, \\
    & I_{it}^{k\omega} \leq C_{it}^{k\omega} & \forall i \in \mathscr{N}, k \in \mathscr{K}, \omega \in \Omega, t \in \mathscr{T} \cup \{T\} \label{invCapacityPH},\\
    & x_{at}^{k\omega} \leq U_a  z_{at}^{k\omega} & \forall a \in \mathscr{A}, k \in \mathscr{K}, \omega \in \Omega, t \in \mathscr{T} \label{arcCapacityPH},\\
    & \sum_{k \in \mathscr{K}} \sum_{t \in \mathscr{T}} y_{at}^{k\omega} \leq 1 & \forall a \in \mathscr{A}, \omega \in \Omega \label{arcConstructionPH} ,\\
    & y_{at}^{k\omega} \leq z_{at}^{k\omega} & \forall a \in \mathscr{A}, k \in \mathscr{K}, \omega \in \Omega, t \in \mathscr{T} \label{arcAssignFirstPH},\\
    & \sum_{k \in \mathscr{K}} \sum_{t' = 0}^t y_{at'}^{k\omega} = \sum_{k \in \mathscr{K}} z_{at}^{k\omega}  & \forall a \in \mathscr{A}, \omega \in \Omega, t \in \mathscr{T} \label{arcAssignLaterPH},\\
    & z_{at}^{k\omega} + \sum_{k' \in \mathscr{K} \setminus \{k\}} z_{a,t-1}^{k'\omega} \leq 1 + b_{at}^{k\omega} & \forall a \in \mathscr{A}, k \in \mathscr{K}, \omega \in \Omega, t \in \mathscr{T} - \{0\} \label{conversionPH},\\ 
    & \sum_{k \in \mathscr{K}} \sum_{t \in \mathscr{T}} b_{at}^{k\omega} \leq 1 & \forall a \in \mathscr{A}, \omega \in \Omega \label{conversionLimitPH}, \\
    & y_{a0}^{k\omega} = 1 & \forall k \in \mathscr{K}, a \in \mathscr{A}_k, \omega \in \Omega \label{initialNetworkPH} ,\\ 
    & I_{i0}^{k\omega} =o_{i}^{k} &  \forall i \in \mathscr{N}, k \in \mathscr{K}, \omega \in \Omega \label{initInvPH} ,\\
    & x_{at}^{k\omega} \geq 0 & \forall a \in \mathscr{A}, k \in \mathscr{K}, \omega \in \Omega, t \in \mathscr{T},\\
    & I_{it}^{k\omega} \geq 0 & \forall i \in \mathscr{N}, k \in \mathscr{K}, \omega \in \Omega, t \in \mathscr{T} \cup \{T\} ,\\
    & e_{it}^{k\omega} \geq 0 & \forall i \in \mathscr{N}, k \in \mathscr{K}, \omega \in \Omega, t \in \mathscr{T},\\
    %& e_{it}^{-k\omega} \geq 0 & \forall i \in \mathscr{N}, k \in \mathscr{K}, \omega \in \Omega, t \in \mathscr{T}\\
    & y_{at}^{k\omega} \in \{0,1\} & \forall a \in \mathscr{A}, k \in \mathscr{K}, \omega \in \Omega, t \in \mathscr{T},\\
    & z_{at}^{k\omega}\in \{0,1\} & \forall a \in \mathscr{A}, k \in \mathscr{K}, \omega \in \Omega, t \in \mathscr{T},\\
    & b_{at}^{k\omega}\in \{0,1\} & \forall a \in \mathscr{A}, k \in \mathscr{K}, \omega \in \Omega, t \in \mathscr{T} \label{phIndvEnd},\\
    & y_{at}^{k\omega} =  y_{at}^{k} & \forall a \in \mathscr{A}, k \in \mathscr{K}, \omega \in \Omega, t \in \mathscr{T} \label{nac1},\\
    & z_{at}^{k\omega} =  z_{at}^{k} & \forall a \in \mathscr{A}, k \in \mathscr{K}, \omega \in \Omega, t \in \mathscr{T} \label{nac2},\\
    & b_{at}^{k\omega} =  b_{at}^{k} & \forall a \in \mathscr{A}, k \in \mathscr{K}, \omega \in \Omega, t \in \mathscr{T}  \label{nac3},\\
    & y_{at}^{k} \in \{0,1\} & \forall a \in \mathscr{A}, k \in \mathscr{K}, \omega \in \Omega, t \in \mathscr{T},\\
    & z_{at}^{k}\in \{0,1\} & \forall a \in \mathscr{A}, k \in \mathscr{K}, \omega \in \Omega, t \in \mathscr{T},\\
    & b_{at}^{k}\in \{0,1\} & \forall a \in \mathscr{A}, k \in \mathscr{K}, \omega \in \Omega, t \in \mathscr{T}, \label{phEnd} 
\end{align}
%Our progressive hedging part of the algorithm is based on \cite{gade2016obtaining}.

As proposed by \cite{rockafellar1991scenarios}, we relax the non-anticipativity constraints by removing Constraints \eqref{nac1} - \eqref{phEnd}, allowing the model for each $\omega \in \Omega$ to be solved individually. \uh{This removes the requirement for first-stage decisions to equal a single global variable (e.g., $y_{at}^{k}$). Instead, the Progressive Hedging algorithm encourages consensus by penalizing each scenario's decisions for deviating from their weighted average. These averages are given below:}
{ \begin{align}
    & \overline{y}_{at}^{k} \coloneqq \sum_{ \omega \in \Omega} p_\omega y_{at}^{k\omega} & \forall a \in \mathscr{A}, k \in \mathscr{K}, t \in \mathscr{T}, \\
    & \overline{b}_{at}^{k} \coloneqq \sum_{ \omega \in \Omega} p_\omega b_{at}^{k\omega} & \forall a \in \mathscr{A}, k \in \mathscr{K}, t \in \mathscr{T}, \\ 
    & \overline{z}_{at}^{k} \coloneqq \sum_{ \omega \in \Omega} p_\omega z_{at}^{k\omega} & \forall a \in \mathscr{A}, k \in \mathscr{K}, t \in \mathscr{T}.
\end{align}}

\uh{Following the framework of \cite{rockafellar1991scenarios}}, The objective for each scenario $\omega$ is then modified with the following Lagrangian coefficients:
\begin{align}
    \text{min  } &  \sum_{a \in \mathscr{A}} \sum_{t \in \mathscr{T}} \sum_{k \in \mathscr{K}} \left( f_{at}^k + \lambda_{at}^{k\omega} \right) y_{at}^{k\omega}  + \left( g_{at}^k + \delta_{at}^{k\omega} \right)b_{at}^{k\omega} + \pi_{at}^{k\omega} z_{at}^{k\omega} \span \nonumber \\ & + \frac{\rho}{2} \left( \left( y_{at}^{k\omega} - \overline{y}_{at}^{k} \right)^2 + \left( b_{at}^{k\omega} - \overline{b}_{at}^{k} \right)^2 + \left( z_{at}^{k\omega} - \overline{z}_{at}^{k} \right)^2 \right)  \nonumber \\  &  +  \sum_{k \in \mathscr{K}} \sum_{t \in \mathscr{T}} \left(  \sum_{a \in \mathscr{A}}   c_{at}^k x_{at}^{k\omega}  + \sum_{i \in \mathscr{N}} h e_{it}^{k\omega} \right). \span \span \label{quadPhObj}
\end{align}
Here, $\lambda$, $\delta$, and $\pi$ are Lagrangian multipliers for decision variables $y$, $b$, and $z$ respectively. The ratio $\rho$ is used to enforce non-anticipativity by penalizing the deviations of decision variables from their average values across scenarios, as well as acting as a step size for the updates on multipliers. %The average values of decision variables over scenarios, $\overline{y}$, $\overline{b}$, and $\overline{z}$, are calculated from the previous iteration of the algorithm. 
Since all the first-stage variables are binary in our model, the quadratic objective function \eqref{quadPhObj} would be equal to the following linearized objective function for all binary solutions:
\begin{align}
    \text{min  } &\sum_{a \in \mathscr{A}} \sum_{t \in \mathscr{T}} \sum_{k \in \mathscr{K}} \left( f_{at}^k + \lambda_{at}^{k\omega} \right) y_{at}^{k\omega}  + \left( g_{at}^k + \delta_{at}^{k\omega} \right)b_{at}^{k\omega} + \pi_{at}^{k\omega} z_{at}^{k\omega} \span \nonumber \\ & + \left( \frac{\rho}{2} - \rho \overline{y}_{at}^{k}\right) y_{at}^{k\omega}  +  \left( \frac{\rho}{2} - \rho \overline{b}_{at}^{k}\right) b_{at}^{k\omega} +  \left( \frac{\rho}{2} - \rho \overline{z}_{at}^{k}\right) z_{at}^{k\omega}\nonumber \\  &  +  \sum_{k \in \mathscr{K}} \sum_{t \in \mathscr{T}} \left(  \sum_{a \in \mathscr{A}}   c_{at}^k x_{at}^{k\omega}  + \sum_{i \in \mathscr{N}} h e_{it}^{k\omega}  \right).\span \span \label{linearizedPhObj}
\end{align}

Making the conversion explicit, $\argmin_{y_{at}^{k\omega}} \frac{\rho}{2} \left( y_{at}^{k\omega} - \overline{y}_{at}^{k} \right)^2 = \argmin_{y_{at}^{k\omega}} \frac{\rho}{2} \left( {y_{at}^{k\omega}}^2 - 2 y_{at}^{k\omega} \overline{y}_{at}^{k} + {\overline{y}_{at}^{k}}^2 \right)$ where ${\overline{y}_{at}^{k}}^2$ might be eliminated as a constant and ${y_{at}^{k\omega}}^2 = y_{at}^{k\omega}$ for $y_{at}^{k\omega} \in \{0,1\}$, resulting in $\argmin_{y_{at}^{k\omega}} \frac{\rho}{2} ( y_{at}^{k\omega} - 2 y_{at}^{k\omega} \overline{y}_{at}^{k} ) = \argmin_{y_{at}^{k\omega}} \left( \frac{\rho}{2} - \rho \overline{y}_{at}^{k}\right) y_{at}^{k\omega}$. The same conversion applies for the rest of the decision variables. Using this reformulation and relaxed model, we describe the steps of the PHA in Algorithm \ref{alg:PH}. %While solving the corresponding models for scenario $\omega$, we assume $\Omega \coloneqq \{\omega\}$ and $p_\omega = 1$.
\begin{algorithm}
    \caption{Progressive Hedging} \label{alg:PH}
    \hspace*{\algorithmicindent} \textbf{Input:} $\lambda \coloneqq 0$, $\delta \coloneqq 0$, $\pi \coloneqq 0$, $\rho > 0$.\\
    \hspace*{\algorithmicindent} \textbf{Output:} A feasible solution to the model \eqref{firstEq} - \eqref{lastEq}. 
    \begin{algorithmic} [1]
        \State\textbf{Initialization:} For each $\omega \in \Omega$, solve \eqref{phStart} - \eqref{phIndvEnd} by assuming $\Omega \coloneqq \{\omega\}$ and $p_\omega = 1$. 
        \While {convergence on first-stage variables is not satisfied}
            \State \textbf{Average:} $\overline{y}_{at}^{k} \coloneqq \sum_{ \omega \in \Omega} p_\omega y_{at}^{k\omega}$, $\overline{b}_{at}^{k} \coloneqq \sum_{ \omega \in \Omega} p_\omega b_{at}^{k\omega}$, $\overline{z}_{at}^{k} \coloneqq \sum_{ \omega \in \Omega} p_\omega z_{at}^{k\omega}$.
            \State \textbf{Price:} $\lambda_{at}^{k\omega} \coloneqq \lambda_{at}^{k\omega} + \rho (y_{at}^{k\omega} - \overline{y}_{at}^{k})$, $\delta_{at}^{k\omega} \coloneqq \delta_{at}^{k\omega} + \rho (b_{at}^{k\omega} - \overline{b}_{at}^{k})$, $\pi_{at}^{k\omega} \coloneqq \pi_{at}^{k\omega} + \rho (z_{at}^{k\omega} - \overline{z}_{at}^{k})$.
            \State \textbf{Solution:} For each $\omega \in \Omega$, solve \eqref{balancePH} - \eqref{phIndvEnd} with the objective \eqref{linearizedPhObj} by assuming $\Omega \coloneqq \{\omega\}$ and $p_\omega = 1$.
        \EndWhile
    \end{algorithmic}
\end{algorithm}

We apply four accelerating techniques to the base PHA in our solution method. The first acceleration technique is to consider bundles of scenarios instead of single-scenario subproblems. This results in solving two-stage recourse models (i.e., subject to the scenarios in the bundle) in the PHA, instead of solving deterministic single scenario problems. It has been shown that for well-calibrated bundle sizes, this reduces computation time without affecting the quality too much \citep{crainic2014scenario}. \uh{For this study, we use random assignment to form scenario bundles, serving as a baseline strategy to test the effectiveness of bundling. %This allows us to isolate the benefit of bundling itself, without the additional complexity of deterministic boundling methods. 
More sophisticated bundling strategies are left for future work, see in Section \ref{sect:concl}.} 
We use a set of bundles $\Xi$ where each scenario $\omega \in \Omega$ is in exactly one bundle $\xi \in \Xi$. The probability of a bundle is assigned as $p_\xi = \sum_{\omega \in \xi} p_\omega$. The objective function of the two-stage recourse model for bundle $\xi$ is given below in Equation \eqref{bundStart}.
\begin{align}
    \text{min  } &  \sum_{a \in \mathscr{A}} \sum_{k \in \mathscr{K}} \sum_{t \in \mathscr{T}} \left( f_{at}^k y_{at}^{k\xi}  + g_{at}^k b_{at}^{k\xi} \right) + \sum_{\omega \in \xi} \frac{p_\omega}{p_\xi} \Bigg[ \sum_{k \in \mathscr{K}} \sum_{t \in \mathscr{T}} \left(  \sum_{a \in \mathscr{A}}   c_{at}^k x_{at}^{k\omega}  + \sum_{i \in \mathscr{N}} h e_{it}^{k\omega} \right) \Bigg] .\label{bundStart}
\end{align}

Secondly, we observe that a majority of first-stage decision variables converge to a common value at early iterations of PHA. In this case, these decision variables may be fixed to their converged value for all scenarios for the remaining iterations of the algorithm. This is called \textit{slamming} in \cite{watson2011progressive} and has shown to be an effective way of accelerating the algorithm without affecting the solution quality substantially. Authors fix variables if they are equal for all scenarios for a given number of consecutive iterations. This given number of iterations is arranged in such a way that it decreases over time to enhance convergence. 

We use a similar acceleration technique with two main differences. Firstly, we enable fixing variables even for the cases where they are not equal for all scenarios but for the majority of scenarios\uh{, i.e., with a soft fixing scheme}. If the total probability of scenarios where the decision variable shares a common value exceeds an input probability parameter $p_{ \textsc{h}}$, we fix the corresponding decision variables. Secondly, we observe that even if some decision variables do not share the same value, we may see similarities between the values of decision variables in terms of our problem context. Thus, we implement a problem-specific convergence rule that does not depend on individual decision variables as in \cite{watson2011progressive}. Instead, the convergence condition depends on subsets of decision variables. For example, consider a case where one arc is built for the majority of scenarios from node $i_1$ to $i_2$ for commodity $k$. \uh{However, since we allow multiple arcs between a node pair $(i_1,i_2)$ ---each representing a pipeline with a different capacity--- arcs with different capacities, or in different periods, may be built in different scenarios.} In this case, the individual values of $y_{at}^{k\omega}$ suggest the algorithm is far from being converged, but the following would be satisfied for the majority of scenarios:
\begin{align}
    \sum_{a \in (i_1,i_2)}  \sum_{t \in \mathscr{T}} y_{at}^{k\omega} \geq 1. \label{heurConst}
\end{align}
This means that for the majority of scenarios, the network for $k$ requires an arc in the flow $(i_1,i_2)$. This is frequently seen, for example, when $(i_1,i_2)$ is one of the few candidate paths connecting some local supply node to a demand point of $k$. 

To implement this logic of convergence, we create the acceleration technique given in Algorithm \ref{alg:HC}. \uh{We refer to this strategy as \textit{aggregate soft fixing}, reflecting both its use of aggregated variable sums and its soft enforcement based on a probability threshold, $p_{ \textsc{h}}$.} The algorithm checks for Constraint \eqref{heurConst} which is defined for decision variable $y$. We observe in preliminary experiments that the complexity of our model is mainly driven by the $y$ variables. As such, there is no need to consider the same logic for the variables $z$ and $b$. Indeed, they can be found rather easily for a fixed set of $y$ variables. This procedure is run once at each iteration of PHA. 

\uh{This approach offers three main advantages. First, it enables earlier fixing than hard fixing techniques, since hard fixing requires $100\%$ agreement across scenarios; in contrast, aggregate soft fixing activates when the total probability of agreeing scenarios exceeds a threshold. Second, by focusing on aggregated subsets rather than individual variables, convergence is detected earlier and more broadly. Third, this aggregation leads to looser constraints—such as $\sum y \geq 1$ instead of strict equalities like $y=1$—giving the solver more flexibility and helping it maintain near-optimality% despite the heuristic nature of the fix
. Finally, while tailored to our network design context, this idea is not problem-specific: many two-stage recourse models in operations research contain similar aggregate patterns. We believe this makes aggregate soft fixing a promising extension for other progressive hedging applications, and we highlight this in Section \ref{sect:concl} as a direction for future work.}

\begin{algorithm}
    \caption{  softFix} \label{alg:HC}
    \hspace*{\algorithmicindent} \textbf{Input:} A probability level of $p_{ \textsc{h}}$.
    \begin{algorithmic} [1]
        \ForEach{pair of nodes $(i_1,i_2)$}
            \ForEach{commodity $k$}
                \If{total probability of bundles satisfying Constraint \eqref{heurConst} is at least $p_{ \textsc{h}}$}
                    \State Add Constraint \eqref{heurConst} to all bundles.
                \EndIf
            \EndFor
        \EndFor
    \end{algorithmic}
\end{algorithm}

The third acceleration technique enables early termination of the PHA. As mentioned above, the majority of first-stage decision variables converge rather fast in the PHA. Often, the algorithm spends the majority of the iterations for some small discrepancies in the last decision variables. In order to accelerate, we terminate the algorithm before it is fully converged, and solve the rest of the problem with a MIP solver. One way of doing this would be fixing a subset of first-stage decision variables, and then solving the MIP for the rest. However, we implement a similar idea to the \uh{aggregate soft fixing} and create heuristic constraints without fixing individual decision variables. It has the same heuristic constraint structure of Algorithm \ref{alg:HC}, but for early convergence, we also check the cases where Constraint \eqref{heurConst} sums to zero. At each iteration, Algorithm \ref{alg:EC} is executed in our method. Once the number of these heuristic constraints reaches the target stopping condition, we solve a MIP for the base model including these constraints. %We again only check on decision variables $y$ as these are deriving the complexity as explained above. %This also resolves the issue that progressive hedging algorithm may not fully converge on some instances, and end up cycling. Progressive hedging converges in linear time on a linear programming models, but same does not apply for mixed-integer programming models.

\begin{algorithm}[h]
    \caption{earlyConvergence} \label{alg:EC}
    \hspace*{\algorithmicindent} \textbf{Input:} A probability level of $p_{ \textsc{e}}$, $i_t \coloneqq 0$, $i_c \coloneqq 0$.\\
    \hspace*{\algorithmicindent} \textbf{Output:} The stopping condition flag (Y/N), and the set of heuristic constraints. 
    \begin{algorithmic} [1]
        \ForEach{pair of nodes $(i_1,i_2)$}
            \ForEach{commodity $k$}
                \State $i_t \coloneqq i_t  + 1$.
                \If{$\sum_{a \in (i_1,i_2)}  \sum_{t \in \mathscr{T}} y_{at}^{kw} = 0$ for all scenarios}
                    \State $i_c \coloneqq i_c  + 1$.
                    \State Record heuristic constraint as $\sum_{a \in (i_1,i_2)}  \sum_{t \in \mathscr{T}} y_{at}^{kw} = 0$.
                \EndIf
                \If{$\sum_{a \in (i_1,i_2)}  \sum_{t \in \mathscr{T}} y_{at}^{kw} \geq 1$ for all scenarios}
                    \State $i_c \coloneqq i_c  + 1$.
                    \State Record heuristic constraint as $\sum_{a \in (i_1,i_2)}  \sum_{t \in \mathscr{T}} y_{at}^{kw} \geq 1$.
                \EndIf
            \EndFor
        \EndFor
        \If{$i_c / i_t \geq p_{\textsc{e}}$}
            \State Return true with the set of heuristic constraints.
        \Else
            \State Return false.
        \EndIf    
    \end{algorithmic}
\end{algorithm}

Finally, since the major complexity of the problem is due to decision variables $y$, as the fourth acceleration technique, we implement the progressive hedging algorithm only on decision variables $y$. This means, we eliminate calculating the average values of $\overline{b}_{at}^{k}$ and $\overline{z}_{at}^{k}$, and remove the corresponding Lagrangian multipliers, $\delta$ and $\pi$. Therefore, the convergence is checked only on decision variables $y$. This would yield the Lagrangian relaxed objective function in Equation \eqref{linearizedBundledObj}.
\begin{align}
    \text{min  } &\sum_{a \in \mathscr{A}} \sum_{t \in \mathscr{T}} \sum_{k \in \mathscr{K}} \left( f_{at}^k + \lambda_{at}^{k\xi} + \frac{\rho}{2} - \rho \overline{y}_{at}^{k} \right) y_{at}^{k\xi}  + g_{at}^k b_{at}^{k\xi} \nonumber \\  &  + \sum_{\omega \in \xi} \frac{p_\omega}{p_\xi} \Bigg[  \sum_{k \in \mathscr{K}} \sum_{t \in \mathscr{T}} \left(  \sum_{a \in \mathscr{A}}   c_{at}^k x_{at}^{k\omega}  + \sum_{i \in \mathscr{N}}   h e_{it}^{k\omega}  \right) \Bigg].\span \span \label{linearizedBundledObj}
\end{align}

The above-mentioned acceleration techniques are implemented in order to create our progressive-hedging-based matheuristic algorithm given in Algorithm \ref{alg:PHBM}.
\begin{algorithm}
    \caption{Progressive Hedging Based Matheuristic} \label{alg:PHBM}
    \hspace*{\algorithmicindent} \textbf{Input:} $\lambda \coloneqq 0$, $\rho > 0$, $p_{ \textsc{h}}$, $p_{ \textsc{e}}$, $\Xi$.\\
    \hspace*{\algorithmicindent} \textbf{Output:} A feasible solution to the model \eqref{firstEq} - \eqref{lastEq}. 
    \begin{algorithmic} [1]
        \State\textbf{Initialization:} For each $\xi \in \Xi$, solve \eqref{balance} - \eqref{lastEq} with the objective \eqref{bundStart} and $\Omega \coloneqq \xi$. 
        \While {stopping condition flag is false}
            \State \textbf{Average:} $\overline{y}_{at}^{k} \coloneqq \sum_{ \xi \in \Xi} p_\xi y_{at}^{k\xi}$.
            \State \textbf{Price:} $\lambda_{at}^{k\xi} \coloneqq \lambda_{at}^{k\xi} + \rho (y_{at}^{k\xi} - \overline{y}_{at}^{k})$.
            \State \uh{softFix}($p_{ \textsc{h}}$)
            \State \textbf{Update:} For each $\xi \in \Xi$, solve \eqref{balance} - \eqref{lastEq} with the objective \eqref{linearizedBundledObj} and $\Omega \coloneqq \xi$.
            \State earlyConvergence($p_{ \textsc{e}}$)
        \EndWhile
        \State \textbf{Solution:} Solve \eqref{firstEq} - \eqref{lastEq} by the set of heuristic constraints.
    \end{algorithmic}
\end{algorithm}

\section{Computational Experiments} \label{sect:compExp}

%This section evaluates the computational efficiency of our progressive-hedging-based matheuristic, focusing on the speed and scalability of the method. %We compare our two-stage stochastic programming solution (SP) with the expected value solution (EV), which is computed by taking the expected values of stochastic parameters and evaluating its solution for each scenario. This help us to measure the Value of Stochastic Solution (VSS), defined as VSS = EV - SP, the value of taking uncertainty into account in the decision-making processes. 
\uh{This section presents a computational study to evaluate our progressive-hedging-based matheuristic (PHM). We assess its performance against two key benchmarks: the direct solution using Gurobi (GRB) and the solution from the deterministic expected value (EV) problem. The comparison with GRB serves to validate our matheuristic's computational efficiency and its ability to find high-quality solutions for large-scale instances where exact methods may struggle. The comparison with the EV solution, in turn, demonstrates the value of the stochastic approach and allows us to analyze the structural differences in the resulting network designs when uncertainty is explicitly considered.}

In Section \ref{sect:bench}, we introduce the benchmark instances and computational parameters. In Section \ref{sect:soln}, we give an overview of the efficiency of our developed methods%, where we also briefly discuss an alternative solution method based on Benders decomposition
. In Section \ref{sect:Bund}, we test the impact of the bundle size on our matheuristic. In Section \ref{sect:var}, we assess the influence of various uncertainty levels on our matheuristic. In Section \ref{sect:timing}, we discuss the difference in solution structure when accounting for uncertainty compared to taking the expected net supply. %\uh{In Section \ref{sect:I}, we examine the impact of omitting the storage option.} 
The case study is presented in Section \ref{sect:Case}.

\subsection{Benchmark Instances and Computational Parameters} \label{sect:bench}
All the experiments are performed using an AMD 7763 CPU (with 16 cores @ 2.45 GHz) with 32GB RAM. All algorithms are executed in C++20 in combination with Gurobi 10.0.1. 
We generate $20$ instances for some combinations of $N \in \{8,10,12,15,20\}$, $T \in \{8,10,15\}$, and $W \in \{60,120,240\}$. \uh{%While no formal in-sample stability analysis was conducted, t
These scenario counts were chosen to balance uncertainty representation and computational tractability. They also align with common practices in stochastic network design, such as \citet{crainic2021partial}, who use $96–256$ scenarios for comparable problem sizes. The values of $W$ are also selected for their high divisibility, enabling testing of various bundle sizes. Based on these tests in Section~\ref{sect:Bund}, we use bundle sizes of $6$ for the case $W=60$, and $8$ for the cases $W=\{120,240\}$.} \uh{We consider a setting with two commodities: one representing natural gas and one representing hydrogen. These are treated as distinct, homogeneous energy types without differentiation by quality grade.}
The number of arcs, $A$, varies per each instance by the instance generator. The average number of arcs for these instances is reported later in our computational results table. 

Each of these instances is solved by our progressive-hedging-based matheuristic (PHM) and Gurobi (GRB). We also obtain the expected value solution (EV) by solving a single-scenario expected value problem and then simulating the performance of this solution over all scenarios in $\Omega$.  For PHM, we apply a MIP gap of $10 \%$ for the \textit{Initialization} and \textit{Update} steps of Algorithm \ref{alg:PHBM}. We further apply a time limit of $1$ hour and a MIP gap of $1 \%$ for the \textit{Solution} step. We set $p_{\textsc{h}}$ and $p_{\textsc{e}}$ to $0.2$ and $1$, respectively. We multiply these probabilities by $0.97$ in each iteration of the algorithm for earlier convergence. We solve the base model \eqref{firstEq}-\eqref{lastEq} via GRB. The EV solution is obtained by taking the average of all stochastic parameters to create an average scenario for each instance. Both GRB and EV have a $3$-hour time limit. The GRB MIP gap is set as $1 \%$. We believe this is a competitive setting to test for GRB, since without these limits it either takes too much time to solve the problem or the chances of going out of memory are too high.

For further details on instance generation and scenario construction, please consult \ref{instGen}.

\subsection{Performance of the Matheuristic} \label{sect:soln}
In this section, we compare our progressive-hedging-based matheuristic (denoted as PHM) with the MIP solution of the base model (denoted as GRB), and the expected value solution (denoted by EV).  \uh{Initially, we explored using Benders decomposition; however, preliminary experiments indicated this approach was not computationally competitive for our problem. Its details are given in \ref{benders}.}
A summary of the average performance over all instances is provided in Table~\ref{table:00}. 

{
\renewcommand{\arraystretch}{0.85}
\begin{table}[h]
\centering
\caption{  Solution quality (Gap \%) and time (Time \%) comparison for the matheuristic (PHM) and Expected Value (EV) solutions. All values are relative to a benchmark solution (either GRB or PHM), with precise definitions provided in the subsequent text.}
\label{table:00}
%\resizebox{\textwidth}{!}{
\begin{tabular}{rrrrrrrrrrrrrr} \toprule
& & & & & & \multicolumn{4}{c}{PHM} & \multicolumn{4}{c}{EV}\\
\cmidrule(r){7-10}\cmidrule(r){11-14}\\
& & & & & & \multicolumn{2}{l}{Gap ($\%$)} & \multicolumn{2}{l}{Time ($\%$)} & \multicolumn{2}{l}{Gap ($\%$)} & \multicolumn{2}{l}{Time ($\%$)}\\
\cmidrule(r){7-8}\cmidrule(r){9-10}\cmidrule(r){11-12}\cmidrule(r){13-14}
$N$ & $A$& $T$ & $W$ & $\#_{\textsc{GRB}}$ & $\#_{\textsc{PHM}}$ & min   & avg  & min  & avg & min   & avg  & min  & avg   \\\midrule
8  & 79  & 8  & 60  & 20 & 20 & -0.31  & 1.34  & 0.56  & 3.21  & 4.28   & 12.35 & 0.13   & 0.72   \\
   &     &    & 120 & 20 & 20 & -0.33  & 1.61  & 0.40  & 1.59  & 4.50   & 12.15 & 0.05   & 0.34   \\
   &     & 10 & 60  & 20 & 20 & 0.16   & 1.24  & 0.60  & 2.32  & 4.27   & 12.32 & 0.09   & 1.19   \\
   &     &    & 120 & 14 & 20 & 0.01   & 1.46  & 0.83  & 3.09  & 4.06   & 10.46 & 0.08   & 0.78   \\
10 & 104 & 8  & 60  & 20 & 20 & 0.70   & 2.52  & 0.44  & 5.65  & 6.17   & 16.91 & 0.02   & 0.93   \\
   &     &    & 120 & 16 & 20 & 0.10   & 1.83  & 0.67  & 5.38  & 5.28   & 16.29 & 0.03   & 0.26   \\
   &     & 10 & 60  & 19 & 20 & 0.51   & 1.70  & 0.64  & 7.11  & 6.67   & 17.36 & 0.05   & 1.39   \\
   &     &    & 120 & 10 & 20 & -10.87 & 0.00  & 2.09  & 8.50  & 5.24   & 15.63 & 0.05   & 0.77   \\
12 & 134 & 8  & 60  & 20 & 20 & 0.33   & 1.41  & 0.45  & 6.72  & 8.00   & 20.23 & 0.05   & 0.76   \\
   &     &    & 120 & 9  & 20 & -0.44  & 1.81  & 1.13  & 14.08 & 6.63   & 18.46 & 0.04   & 1.14   \\
   &     & 10 & 60  & 20 & 20 & 0.35   & 1.91  & 0.72  & 11.14 & 10.13  & 21.16 & 0.19   & 2.03   \\
   &     &    & 120 & 10 & 20 & -8.49  & 0.27  & 1.82  & 23.44 & 6.70   & 17.42 & 0.16   & 8.10   \\
15 & 179 & 10 & 120 & 16 & 19 & -17.06 & -8.98 & 9.22  & 28.94 & 8.19   & 18.26 & 1.39   & 31.58  \\
   &     &    & 240 & 0  & 20 & -      & -     & 13.69 & 36.88 & -53.25 & 11.57 & 0.80   & 27.70  \\
20 & 258 & 10 & 240 & 0  & 6  & -      & -     & 10.78 & 35.43 & 6.17   & 9.26  & 1.09   & 71.95  \\
   &     & 15 & 240 & 0  & 1  & -      & -     & 47.50 & 47.50 & 10.38  & 10.38 & 1.63   & 91.73       
\\ \bottomrule                  
\end{tabular}
%}
\end{table}
}

Each row in Table \ref{table:00} denotes the performance over 20 instances. Column $A$ denotes the average number of arcs over the $20$ instances. The columns $\#_{\textsc{GRB}}$ and $\#_{\textsc{PHM}}$ are the number of instances that are solved without any memory issues. We do not report the number of solutions for the EV, since it stays within the memory limit for all instances. \uh{The remaining columns compare PHM and EV with GRB in terms of objective value (Gap \%) and solution time (Time \%). Gaps are computed using the best-found feasible solutions, using only the instances solved by both methods without memory issues. For PHM, the gap is the percentage difference from the GRB solution. For EV, the gap is the percentage difference to the best solution found by either GRB or PHM; thus in cases where GRB fails on all instances, the PHM solution fully serves as the gap benchmark. ``Time (\%)" columns report the solution time relative to the GRB run time. In cases where $\#_{\textsc{GRB}} = 0$, the time is instead compared against the 3-hour time limit.}

We observe that our algorithm can find on average solutions that are within $1-2\%$ of optimality by only spending $5-6\%$ of the MIP solution time for smaller instances. Moreover, our algorithm outperforms GRB as shown by the negative value for ``Gap" for some of the instances. This negative gap occurs because GRB operates under a time limit, preventing unrealistic computational times and memory issues. Still, our algorithm has fewer memory issues. This is significant for larger-sized instances, as most of these instances cannot be solved by GRB efficiently. %The time ratio on average increases with larger-sized instances, but this is mostly because GRB has a fixed $ 3$-hour time limit. 
This results in an average $8.98\%$ improvement by PHM compared to GRB for the instances with $N=15$, $T=10$, and $W=120$, while PHM spends on average $28.94\%$ of the time GRB takes on these instances. 

We also make comparisons with the EV. In the column ``Gap" of EV, the EV solution is compared with the best solution found via GRB and PHM. Note for the cases that GRB cannot solve any instance, this is always the solution provided by PHM. The solution time of EV is the percentage of GRB time, except for the cases $\#_{\textsc{GRB}} = 0$ where we compare the EV time by $3$ hours. While EV works fast on small instances, we observe an overall gap of $10-20\%$ in its solution quality. Moreover, there is an increase in its time compared to GRB on larger instances, since GRB starts hitting the $3$-hour limit more and more. PHM has similar time ratios on the same larger instances while having better results on the gap. In order to better understand the effect of computational limits on performance comparisons, we conducted additional tests on instances with $N=10$, $T=8$, and $W=60$, but without any time limits and MIP gaps, and using 70GB of RAM for both GRB and PHM. For these specific instances, both the gap and the time required decreased for PHM on average compared to GRB; the gap reduced from $2.52\%$ to $1.46\%$, and the time ratio from $5.65\%$ to $1.50\%$. These results support that imposing limits is beneficial for maintaining GRB’s efficiency compared to PHM.

%On the largest instances with $N=25$, $T=20$, and $W=240$, we observe that all $20$ instances hit the $3$-hour time limit even with $1$-scenario EV solution, due to the immense size of the problem.

It is important to note that for the instances with $N=20$, PHM improves the EV solution by about $10\%$ while spending half of the EV solution time. This is remarkable since while the EV is a $1$-scenario problem, the PHM deals with thirty $8$-scenario problems on $W=240$. The algorithm's success can be attributed not only to the high MIP gap but also to the \uh{aggregate soft fixing} through heuristic constraints, which enhances efficiency. Each iteration is solved as fast as possible, and at the end of each iteration, we fix decision variables as much as possible while keeping overall solution quality high. Only in one extreme case on $N=15$, $T=10$, and $W=240$, we observe that the EV is significantly better (due to PHM not converging quickly enough) as it provides a solution with a  $-53.25\%$ gap. However, even when this is included, on average these instances are improved by $11.57\%$ in terms of the solution quality by PHM.

\subsection{Impact of Scenario Count per Bundle on Solution} \label{sect:Bund}

\begin{figure}[h]
     \centering
     \begin{subfigure}[t]{0.48\textwidth}
        \centering
        \resizebox{0.9\textwidth}{!}{%
        \includegraphics{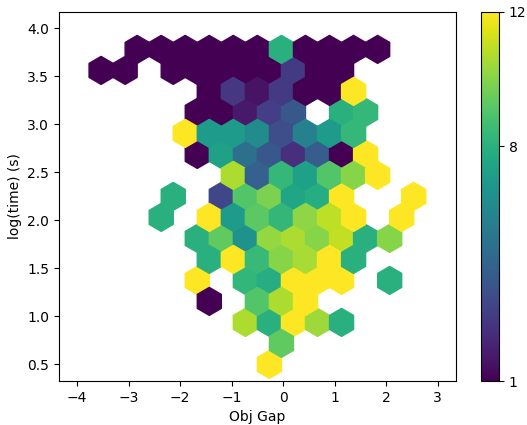}
        }
        \caption{  Hexbin map of individual solutions for 200 instances with $(N,T,W) = (10,8,120)$ under different bundling sizes ($|\xi| \in {1, 8, 12}$). Color intensity corresponds to the concentration of solutions.}
        \label{fig:bundling1}
     \end{subfigure}
     \hfill
     \begin{subfigure}[t]{0.48\textwidth}
         \centering
         \begin{tikzpicture}[scale=0.8]
        \begin{axis}[
            xlabel={$|\xi|$},
            ylabel={Obj},
            xtick=data,
            xtick pos=bottom, 
            yticklabel style={/pgf/number format/fixed, text=red},
            yticklabel pos=left, % Position y-axis labels on the left side
            ytick align=inside, % Align y-axis labels inside the plot area
            ytick pos=left, % Position y-axis ticks on the left side
            ylabel style={text=red}, % Set y-label text color to red
            y axis line style={red}, % Set y-axis line color to red
        ]
        \addplot[color=red,mark=x, smooth] coordinates {
		(1,771.26)
(2,772.9)
(3,772.66)
(4,772.36)
(6,772.96)
(8,773.65)
(12,775)
(24,777.27)
	};
        %\legend{Obj}
        \end{axis}
        
        \begin{axis}[
            axis y line*=right,
            ylabel={Time (sec)},
            ylabel near ticks,
            ylabel style={yshift=-10pt}, 
            xtick=\empty, 
            ticklabel pos=right, % Position y-axis labels on the right side
            ytick align=inside, % Align y-axis labels inside the plot area
            ytick pos=right, % Position y-axis ticks on the right side
            ylabel style={text=blue}, % Set y-label text color to red
            y axis line style={blue}, % Set y-axis line color to red
            yticklabel style={/pgf/number format/fixed, text=blue},
        ]
        \pgfplotsset{every axis y label/.append style={rotate=180, yshift=10pt}}
        \addplot[color=blue,mark=*, smooth] coordinates {
		(1,399.73)
(2,198.39)
(3,131.75)
(4,104.76)
(6,70.31)
(8,68.61)
(12,58.06)
(24,49.98)
};
        %\legend{Time}
        \end{axis}
    \end{tikzpicture}
    \caption{  Average objective (left axis) and solution time (right axis) over a combined set of 200 instances (100 of size $(N,T,W) = (5,8,120)$ and 100 of size $(N,T,W) = (8,8,120)$) for varying bundle sizes $|\xi| \in \{1,2,3,4,6,8,12,24\}$}
    \label{fig:bundling2}
     \end{subfigure} 
        \caption{  Impact of the number of scenarios ($|\xi|$) on solution quality and computation time across different instance sets.}
        \label{fig:bundling}
\end{figure}

%The existing literature on progressive hedging methods discusses how scenario bundling affects solution quality \citep{crainic2014scenario}. 
This subsection analyzes how varying the number of scenarios per bundle affects the solution quality of the PHM. We randomize various instances by using the same parameter settings as before.  We test various numbers of scenarios per bundle, i.e. $ |\xi| \in \{1,2,3,4,6,8,12,24\}$ to randomly assign the scenarios into bundles\uh{, once for each combination of instance and $|\xi|$}. We report their solution in terms of the quality and the solution time in Figure \ref{fig:bundling}, where Figure \ref{fig:bundling1} and Figure \ref{fig:bundling2} provide the individual and the average behavior, respectively.

We provide a hexagon scatter plot in Figure \ref{fig:bundling1} for a subset of the number of scenarios considered. We solve a total of $200$ instances with $(N,T,W) = (10,8,120)$ for each value $ |\xi|  \in \{1,8,12\}$ by grouping scenarios randomly. The colors represent the composition of data included in the hexagon. The objective gap is normalized; it is the gap to the average objective taken from the same instances for all number of scenarios in $ \{1,2,3,4,6,8,12,24\}$. The time is reported in a logarithmic scale, for example, a score of $1.5$ means it takes $10^{1.5}$ seconds to solve. %the problem under the selected $ |\xi| $ level. 
The figure shows the trade-off between the solution quality and the solution time. While increasing the number of scenarios per bundle decreases the solution quality in terms of the objective value, it also decreases the solution time on average. However, we like to stress that the differences in objective value are small for any bundling.

A similar behavior is seen in Figure \ref{fig:bundling2}. The figure provides the average objective value and corresponding solution times of a total of $200$ instances with either $(N,T,W) = (5,8,120)$ or $(N,T,W) = (8,8,120)$. Each of these instances is solved for various randomized bundles with $ |\xi| \in \{1,2,3,4,6,8,12,\allowbreak 24\} $. We observe that the trade-off observed in Figure \ref{fig:bundling1} persists, except for $ |\xi| \in \{2,3\} $. This trade-off is due to using our soft fixing approach, which incorporates heuristic constraints when a certain percentage of instances share common pipeline arcs. As $|\xi|$ increases, even a few bundles sharing the same arc can trigger the heuristic constraints. In the highest value of $|\xi|=24$, heuristic constraints are introduced even if only one bundle contains a particular arc. Because when $W=120$, $p_{\textsc{h}}= 20\%$ of the solutions consist of exactly one bundle with $|\xi|=24$ scenarios. Thus, with higher values of $|\xi|$, the incorporation of heuristic constraints becomes more frequent. While this restriction on the solution space indeed tends to increase the objective value, it also leads to a decrease in the solution time. 

For $ |\xi| = \{1,4,6,8,12,24\} $, the decision-maker may select the trade-off between solution quality and solution time. Regardless of the chosen trade-off, our method consistently exhibits small objective gaps across all values of $ |\xi|$, underscoring its robustness and efficacy. \uh{To justify our focus on average values, we examined the distribution of objective values across the 200 randomized runs. The results show limited variability: the interquartile range remains below $0.75\%$ of the average objective for all bundle sizes ---except $|\xi| = 24$ which showed a slightly increased range of $1.36\%$--- indicating that random scenario grouping introduces no significant variance. %This supports the robustness of our conclusions based on average performance.
}

\uh{As introduced at the beginning of this section, we select either $ |\xi| = 6 $  or $ |\xi| = 8 $ for our computational experiments}. This choice is influenced by the observation that the average objective gap appears to be relatively insensitive to the number of scenarios in our particular setting, and for $6$ or more values the solution time decreases sharply. However, the objective gap range highly depends on the instance parameters selected. For instance, scenarios characterized by higher uncertainty in supply/demand and high penalty costs may exhibit more substantial differences in objective value. In this case, it might be sensible to decrease $ |\xi| $ to find a better balance between the objective function value and the solution time. Similarly, for larger instances requiring extended computational time, it may be advantageous to increase the number of scenarios $|\xi|$ even further, aiming to find some feasible solutions in reasonable computational time. This is particularly important because as computational time increases, the likelihood of encountering memory issues also rises.

\subsection{Impact of the Uncertainty on Solution Quality}\label{sect:var}

In this section, we examine the impact of the uncertainty level on the performance of our progressive-hedging-based matheuristic (PHM) compared to the Gurobi (GRB). We analyze instances with $N=12$ and $A=134$, on low, normal, and high uncertainty in scenarios to explore their influence on gap and solution time. These solutions are executed in C++20 in combination with Gurobi 11.0.0. Normal uncertainty instances are the same as those in Section \ref{sect:soln}. Low uncertainty instances are created similar to the normal ones, but we remove the uncertainty in supply and intermediate locations for hydrogen (see \ref{instGen} for details). High uncertainty instances are introduced by adjusting the range of variability parameters as described in \ref{instGen}: $[G_L,G_U] = [1.2,1.3]$, $[H_L,H_U] = [1.75,2.25]$, and setting $r$ to be uniformly randomized within $[1.01,1.1]$ for Eq. \eqref{a1}. The results are presented in Table \ref{table:var}.

{
\renewcommand{\arraystretch}{0.85}
\begin{table}[h]
\centering
\caption{  Performance of PHM and EV relative to GRB for instances with $N=12$ and $A=134$, analyzed across Low, Normal, and High uncertainty (Var) settings.}
\label{table:var}
%\resizebox{\textwidth}{!}{
\begin{tabular}{lrrrrrrrrrrrr} \toprule
& & & & & \multicolumn{4}{c}{PHM} & \multicolumn{4}{c}{EV}\\
\cmidrule(r){6-9}\cmidrule(r){10-13}\\
& & & & & \multicolumn{2}{l}{Gap ($\%$)} & \multicolumn{2}{l}{Time ($\%$)} & \multicolumn{2}{l}{Gap ($\%$)} & \multicolumn{2}{l}{Time ($\%$)}\\
\cmidrule(r){6-7}\cmidrule(r){8-9}\cmidrule(r){10-11}\cmidrule(r){12-13}
Var & $T$ & $W$ & $\#_{\textsc{GRB}}$ & $\#_{\textsc{PHM}}$ & min   & avg  & min  & avg & min   & avg  & min  & avg   \\\midrule
Low    & 8  & 60  & 20 & 20 & 0.17   & 1.10  & 0.56  & 2.07  & 4.69  & 17.41 & 0.04 & 0.85  \\
       &    & 120 & 11 & 20 & 0.39   & 1.24  & 1.08  & 6.11  & 5.14  & 15.77 & 0.05 & 1.05  \\
       & 10 & 60  & 19 & 20 & 0.30   & 1.71  & 0.60  & 4.30  & 6.36  & 17.60 & 0.11 & 1.85  \\
       &    & 120 & 10 & 20 & -12.74 & -1.05 & 2.02  & 16.25 & 5.54  & 13.80 & 0.30 & 3.50  \\\midrule
Normal & 8  & 60  & 20 & 20 & 0.03   & 1.41  & 0.44  & 3.72  & 7.90  & 17.46 & 0.04 & 0.92  \\
       &    & 120 & 10 & 20 & 0.33   & 1.36  & 0.75  & 12.24 & 6.74  & 17.01 & 0.07 & 1.63  \\
       & 10 & 60  & 20 & 20 & 0.41   & 1.50  & 0.55  & 14.04 & 10.21 & 18.66 & 0.23 & 3.23  \\
       &    & 120 & 10 & 20 & -12.19 & -0.70 & 1.97  & 19.32 & 9.24  & 15.64 & 0.21 & 6.87  \\\midrule
High   & 8  & 60  & 20 & 20 & 0.12   & 0.86  & 4.90  & 15.39 & 9.57  & 21.54 & 0.19 & 3.55  \\
       &    & 120 & 12 & 20 & -0.53  & 0.52  & 8.89  & 20.74 & 5.81  & 20.15 & 0.72 & 5.23  \\
       & 10 & 60  & 19 & 20 & -0.37  & 0.59  & 13.71 & 33.89 & 8.63  & 20.79 & 3.18 & 36.44 \\
       &    & 120 & 16 & 20 & -1.61  & 0.18  & 15.83 & 36.33 & 8.24  & 19.38 & 1.95 & 54.44     
\\ \bottomrule                  
\end{tabular}

%}
\end{table}
}

The column Var represent the variability, and the remaining columns are the same as those in Section \ref{sect:soln}. With more uncertainty, solution times for all methods increase. GRB frequently reaches its 3-hour time limit, which negatively impacts its solution quality. This leads to an increased time ratio for PHM, but GRB is unable to achieve near-optimal solutions within the set time frame. Therefore, PHM demonstrates lower gaps in solution on average as uncertainty increases (except one extreme case on $T=10$ and $W=120$ resulting a high negative gap for low and normal uncertainty levels). Overall, an increase in uncertainty has effects similar to those observed when increasing the problem size, as discussed in Section \ref{sect:soln}; resulting in lower solution gaps and higher time ratios between GRB and PHM. This emphasizes the robustness of PHM under different uncertainty levels.

%Additionally, we observe some cost implications due to changes in uncertainty levels; specifically, total pipeline construction costs decrease by $23\%$ when uncertainty is lowered compared to the normal case, while high uncertainty triples construction costs, emphasizing the sensitivity of cost outcomes to uncertainty parameters. This high increase is partly due to increased demand expectations, which result from high uncertainty in economic growth projections, requiring more pipeline construction.

\subsection{Timing of the Transition in the Hydrogen Economy} \label{sect:timing}

In this section, we investigate the dynamics of transitioning from a natural gas network to one primarily based on hydrogen. We conduct \uh{experiments using} 100 instances on each of the parameter sets $(N,T,W)=(8,10,60)$, representing a smaller network with fewer nodes and potential connections, and $(N,T,W)=(12,10,60)$, illustrating a larger, more complex network structure. While our instances are designed to evaluate computational performance, they are derived to reflect real-world conditions and closely align with the parameters of actual case studies in the hydrogen economy, see \ref{instGen} for details, making them practical for our analysis in this section. We compare the average proportion of network arcs allocated to hydrogen transportation in each period. All these instances are solved via PHM and EV%, with $|\xi|=6$ for PHM
. Figure \ref{fig:hydrogen-transitionV2} illustrates the evolution of hydrogen implementation over time, where the y-axis represents the fraction of arcs being designated for hydrogen transport.

\begin{figure}[h]
\centering
\begin{tikzpicture}[scale= 0.9]
\begin{axis}[
xlabel=Time Period,
ylabel=Hydrogen Ratio,
ymin=0,
xmin=1,
legend pos=south east,
%ymax=1
]
\addplot[color=blue,mark=*, smooth] coordinates {
(1,0)
(2,0.010369)
(3,0.014488)
(4,0.077752)
(5,0.258842)
(6,0.541305)
(7,0.591366)
(8,0.62867)
(9,0.645662)
(10,0.654279)
};
\addlegendentry{$N=8$, PHM}
\addplot[color=cyan,mark=*, smooth] coordinates {
(1,0)
(2,0.048831)
(3,0.191662)
(4,0.348827)
(5,0.399953)
(6,0.456158)
(7,0.494954)
(8,0.534231)
(9,0.560708)
(10,0.569316)
};
\addlegendentry{$N=8$, EV}
\addplot[color=red,mark=x, smooth] coordinates {
(1,0.002222)
(2,0.012745)
(3,0.030258)
(4,0.143667)
(5,0.402821)
(6,0.609418)
(7,0.645589)
(8,0.672319)
(9,0.691441)
(10,0.704949)
};
\addlegendentry{$N=12$, PHM}
\addplot[color=BurntOrange,mark=x, smooth] coordinates {
(1,0)
(2,0.068388)
(3,0.259609)
(4,0.40342)
(5,0.452792)
(6,0.501203)
(7,0.539116)
(8,0.572164)
(9,0.602051)
(10,0.609412)
};
\addlegendentry{$N=12$, EV}
’\end{axis}
\end{tikzpicture}
\caption{  Evolution of the hydrogen network over time, comparing the proportion of network arcs designated for hydrogen for PHM and EV solutions across two network sizes (N=8 and N=12)}
\label{fig:hydrogen-transitionV2}
\end{figure}

As depicted in Figure \ref{fig:hydrogen-transitionV2}, by the tenth period, about $65-70\%$ of network arcs should be designated for hydrogen transport according to the PHM solution. %This is roughly in line with the envisioned $70-75\%$ retrofitted pipelines of the European Hydrogen Backbone, anticipated by the time the full hydrogen transition is completed around 2050 \citep{wang2020european}.
Nevertheless, roughly $30\%$ of the natural gas infrastructure endures even though we assume that the natural gas economy mostly phases out by $t=10$ in the instances. This is due to the inherent incompatibilities between certain natural gas arcs and the evolving hydrogen transport routes. These incompatibilities arise from various factors. Firstly, some natural gas pipelines may require new counterparts to serve new demands arising in the future hydrogen economy, for example with pipelines at a different capacity. Secondly, certain natural gas supplier locations at $t=0$ may not serve as hydrogen suppliers in the future. In such cases, the pipelines originating from these locations are less likely to be converted for hydrogen transport, leading to the existence of pipelines that are no longer in use at $t=10$. This also partly explains that the end ratio of hydrogen pipelines in the case of $N=12$ is about $5\%$ higher than in the case of $N=8$. With more locations in the problem instance, there are increased opportunities for retrofitting existing natural gas infrastructure for hydrogen transport due to the denser network with more supply and demand locations, resulting in a more effective transition. %We discuss this further in Section \ref{sect:Case}.

It is also noteworthy that there is a difference in the average behaviors between the two solution approaches, PHM and EV. The EV suggests converting pipelines earlier in the transition process, with approximately $40\%$ of arcs allocated to hydrogen by $t=4$, whereas PHM only has about $10\%$ at the same time in our instances. Conversely, the EV solutions result in about $8\%$ less conversion compared to the PHM at $t=10$. These disparities arise primarily due to the way the uncertainty is handled. PHM adopts a multi-scenario solution approach, which inherently accounts for a range of potential scenarios, including extreme cases. In contrast, EV employs a deterministic equivalent of the problem, assuming a single-scenario (mean) perspective, which leads to a more optimistic outlook. With this mean-value perspective, EV solutions tend to initiate pipeline construction earlier, assuming a degree of certainty in decision-making. Meanwhile, PHM focuses on optimizing the timing of pipeline construction to minimize overall costs, often favoring delayed pipeline construction ---which is cheaper due to the assumption that costs decrease over time--- while paying the penalty of unsatisfied demand for the extreme cases at early periods. A similar pattern is seen, for example, in other investment decisions under uncertainty. Deterministic approaches might lean towards early investments due to anticipated benefits, while stochastic methods delay them to hedge against potential risks \citep[see, e.g.,][]{dixit1994investment}. 

Nevertheless, due to the additional uncertainty information considered, PHM's solution proves to be more cost-efficient in the long run. Figure \ref{fig:54-1} and \ref{fig:54-2} reflect the average per period cost among two sets of instances, $(N,T,W)=(8,10,60)$ and $(N,T,W)=(12,10,60)$, where EV values are from simulations with suggested construction policies. The costs are provided as ratios to the average cost of PHM and EV at $t=1$. Figure \ref{fig:54-1} demonstrates that PHM increasingly succeeds in satisfying the demand through flow costs as the planning horizon progresses. The widening gap in flow costs indicates a better covering of the demand of the network, which results in PHM incurring lower penalty costs compared to EV, especially evident in later periods as shown in Figure \ref{fig:54-2}. We also observe a periodic pattern in costs due to seasonal factors, which is further explained in \ref{instGen}. The observed differences are most pronounced towards the end of the time horizon, as period $t=10$ mimics sustained long-term use beyond the planning horizon, with increased penalty costs for any unmet demand at $t=10$.

\begin{figure}[t]
    \centering
    \begin{subfigure}[b]{0.49\textwidth}
        \centering

        \begin{tikzpicture}[scale= 0.9]
\begin{axis}[
xlabel=Time Period,
ylabel=Flow Cost Ratio,
ymin=0,
xmin=1,
legend pos=south east,
%ymax=1
]
\addplot[color=blue,mark=*, smooth] coordinates {
(1,1.02571)
(2,1.04176)
(3,0.947138)
(4,1.029533)
(5,0.988384)
(6,1.288805)
(7,1.336172)
(8,1.697992)
(9,1.8442)
};
\addlegendentry{PHM}
\addplot[color=BurntOrange,mark=*, smooth] coordinates {
(1,0.97429)
(2,1.019416)
(3,0.938459)
(4,1.003768)
(5,0.859079)
(6,0.923968)
(7,0.749908)
(8,0.847111)
(9,0.786324)
};
\addlegendentry{EV}
\end{axis}
\end{tikzpicture}

        \caption{  Flow cost ratio over time.}
        \label{fig:54-1}
    \end{subfigure}
    \hfill
    \begin{subfigure}[b]{0.49\textwidth}
        \centering

        \begin{tikzpicture}[scale= 0.9]
\begin{axis}[
xlabel=Time Period,
ylabel=Penalty Cost Ratio,
ymin=0,
xmin=1,
legend pos=south east,
%ymax=1
]
\addplot[color=blue,mark=*, smooth] coordinates {
(1,1.148119)
(2,1.23559)
(3,1.132963)
(4,1.742112)
(5,1.570174)
(6,2.365162)
(7,1.879508)
(8,2.762795)
(9,1.702862)
};
\addlegendentry{PHM}
\addplot[color=BurntOrange,mark=*, smooth] coordinates {
(1,0.851881)
(2,1.405499)
(3,1.151671)
(4,1.952717)
(5,2.117617)
(6,3.442406)
(7,3.32039)
(8,5.108038)
(9,4.247076)
};
\addlegendentry{EV}
’\end{axis}
\end{tikzpicture}

        \caption{  Penalty cost ratio over time.}
        \label{fig:54-2}
    \end{subfigure} 
    \caption{  Average per-period costs for PHM and EV solutions, shown as a ratio relative to the average cost at t=0.}
    \label{fig:54}
\end{figure}

In summary, the interplay between infrastructure expansion and emerging supply demands adds depth to our understanding of the gradual transition from natural gas to hydrogen within the network. %Our analysis of the transition dynamics from a natural gas network to a hydrogen-based infrastructure yields several noteworthy insights. During our analysis, 
We specifically observe a series of periods, mainly occurring in the middle of the planning horizon, marked by a rapid and pronounced transition. Planning and budgeting for these middle periods will be crucial, as they require more resources. This aligns with current energy transition plans in the HEAVENN project as well, where a small adoption of hydrogen until 2040, followed by a rapid growth phase from 2040 to 2050 is expected. Our results mirror this trend.

\section{Case Study and Discussion: Transitioning Natural Gas Network to Hydrogen in the HEAVENN Region} \label{sect:Case}

\begin{figure}[h]
    \centering
    \includegraphics[scale=0.5]{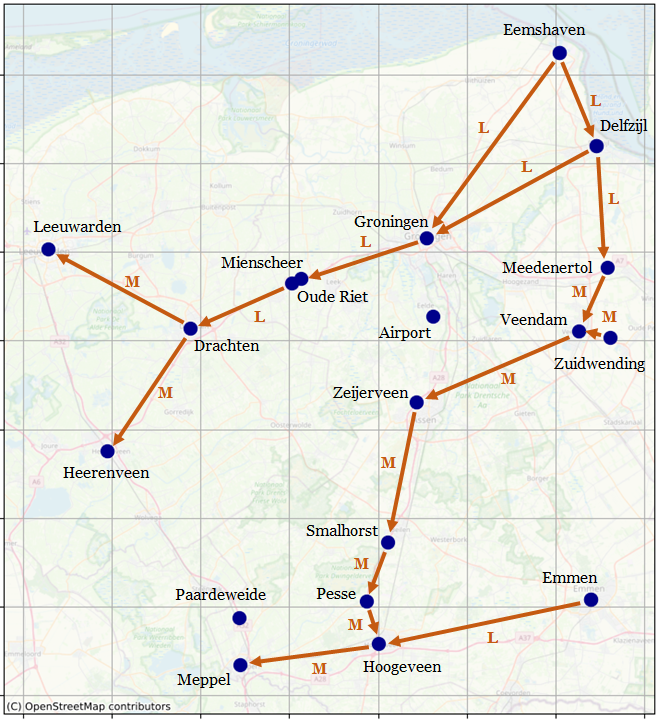}
    \caption{  The initial configuration of the natural gas network for the HEAVENN case study. The letters on the pipeline arcs represent their capacity: S for Small, M for Medium, and L for Large.}
    \label{fig:aprime}
\end{figure}

In this section, we discuss the implications of transitioning the extensive natural gas network within the Northern Netherlands region into a green hydrogen distribution infrastructure. This is a critical step towards achieving a climate-neutral economy as is also a part of the project \textit{Hydrogen Energy Applications in Valley Environments for Northern Netherlands} \citep{heavenncite}, which is Europe's first hydrogen valley initiative. While new hydrogen pipelines are being developed, there is a parallel effort to repurpose the existing natural gas network for hydrogen transport, a cost-effective and environmentally friendly alternative. Based on discussions with our project partners, we focus on a set of locations representing the main critical areas in the network, as shown in Figure \ref{fig:aprime}. These locations encompass five distinct types: suppliers, hydrogen refueling stations or analogous demand units for natural gas, residential areas with heating gas demand, industrial zones requiring gas as a feedstock, and intermediate/storage points facilitating pipeline connections. 

We investigate two variants of the network structure to capture different supply strategies that may emerge in the region. The first variant relies on a single supplier, the Eemshaven port in our specific case, representing the import of hydrogen from abroad, a possible scenario foreseen for regions without sufficient local hydrogen production facilities. In the second variant, we study a situation where several smaller suppliers are spread across the Northern Netherlands region, producing hydrogen locally. For clarity, we refer to the first variant as \textit{Eemshaven import}, and the second variant as  \textit{decentralized production}. While we anticipate a blend of both import and local production, our focus on these two settings aims to offer distinct insights into their potential consequences. The \textit{Eemshaven import} variant is characterized by the certainty at the supply location, unlike the \textit{decentralized production} variant which has variability in supplier locations across scenarios. Both variants start with the initial natural gas network structure represented in Figure \ref{fig:aprime}. On this starting network, we consider a subset of major pipelines in the region. %We select a subset of major pipelines connecting these locations and likely to yield meaningful insights for our research. 
Capacities for these pipelines are categorized into three levels: S, M, and L, representing small, medium, and large capacities, respectively.

We assign parameter values based on consultations with our industry partners in the Northern Netherlands region and from prior studies conducted within the HEAVENN project \citep{robertMSC, julianMSC, hasturk2024stochastic}. Key parameters include information on pipeline transition costs, possible pipeline capacities, initial network of natural gas, seasonality on demand, and the expected pace of the transition from natural gas to hydrogen. The majority of the parameters align with those used in the computational experiments (see also \ref{instGen}). Distinctively, we adjust pipeline conversion cost as $30\%$ of the fixed cost of construction of the same pipeline. % from scratch in our case study to investigate its effects on a real case. 
For the progressive hedging algorithm, we employ the same parameters as in the computational experiments. Our instance generation method is detailed in \ref{instGen}, and it is implemented with $60$ scenarios. 

\subsection{Results} \label{sec:1supp}
In our analysis, we address two distinct network supply settings. As previously outlined, the \textit{Eems\-haven import} variant involves a single supplier, situated at the Eemshaven, representing the import of hydrogen from external sources. In contrast, the \textit{decentralized production} variant explores a scenario where multiple local suppliers are dispersed across the Northern Netherlands region, with a specific focus on uncertainty regarding their locations. The candidate locations for these suppliers are Eemshaven, Delfzijl, Groningen, Zeijerveen, Zuidwending, and Eelde Airport. We randomly assign these locations as supply nodes or intermediary nodes in each of the $60$ scenarios. Additionally, this allocation process also establishes whether these nodes can function as storage facilities. 

\begin{figure}[h]
    \centering
    \begin{subfigure}[b]{0.49\textwidth}
        \centering
        \includegraphics[scale=0.38]{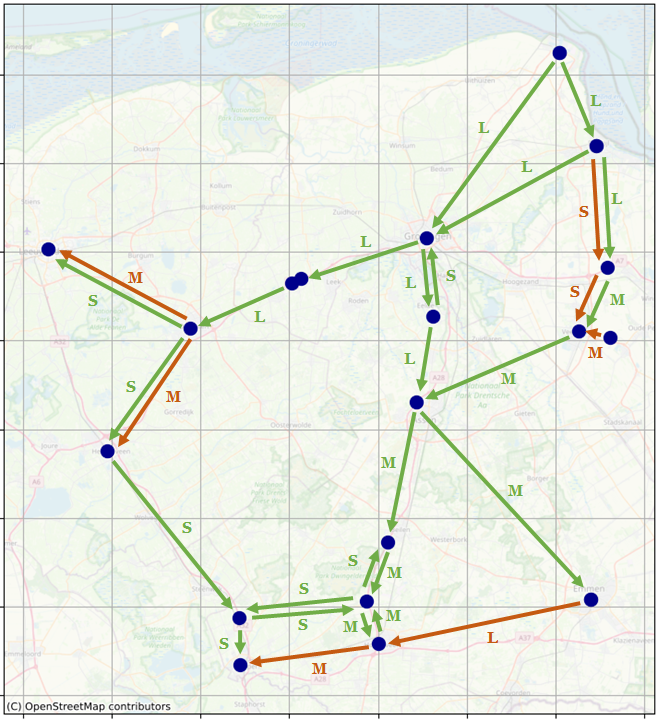}
        \caption{Eemshaven import}
        \label{fig:case1}
    \end{subfigure}
    \hfill
    \begin{subfigure}[b]{0.49\textwidth}
        \centering
        \includegraphics[scale=0.38]{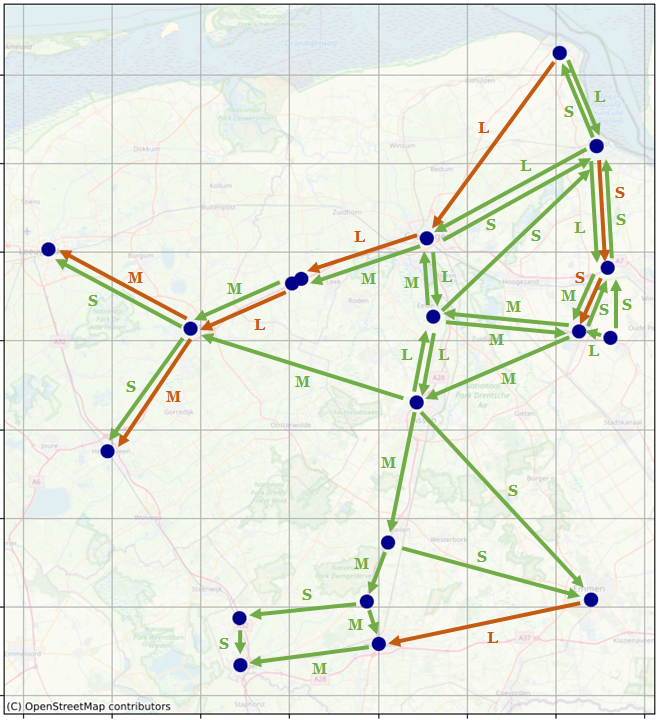}
        \caption{Decentralized production}
        \label{fig:case2}
    \end{subfigure}
    \vfill
    \begin{subfigure}[b]{0.49\textwidth}
        \centering
        \includegraphics[scale=0.38]{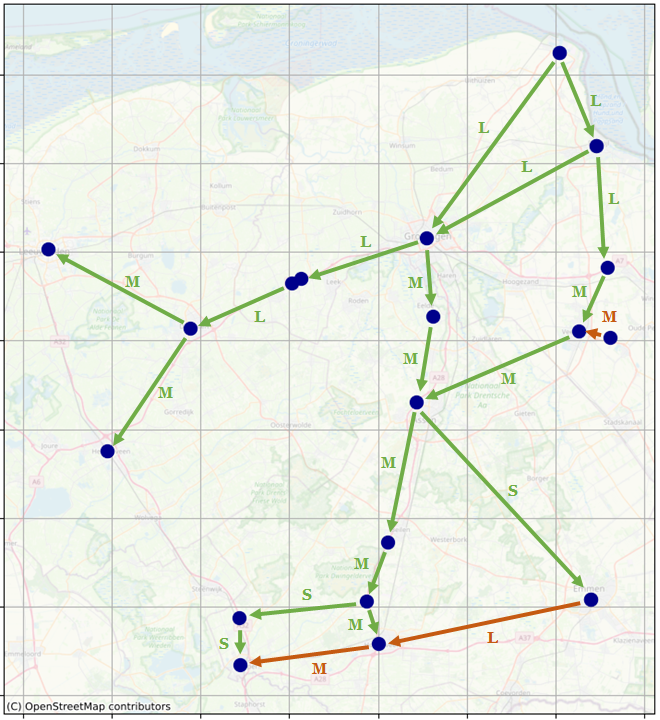}
        \caption{Eemshaven import on deterministic case}
        \label{fig:case1-eev}
        \end{subfigure}
    \hfill
    \begin{subfigure}[b]{0.49\textwidth}
        \centering
        \includegraphics[scale=0.38]{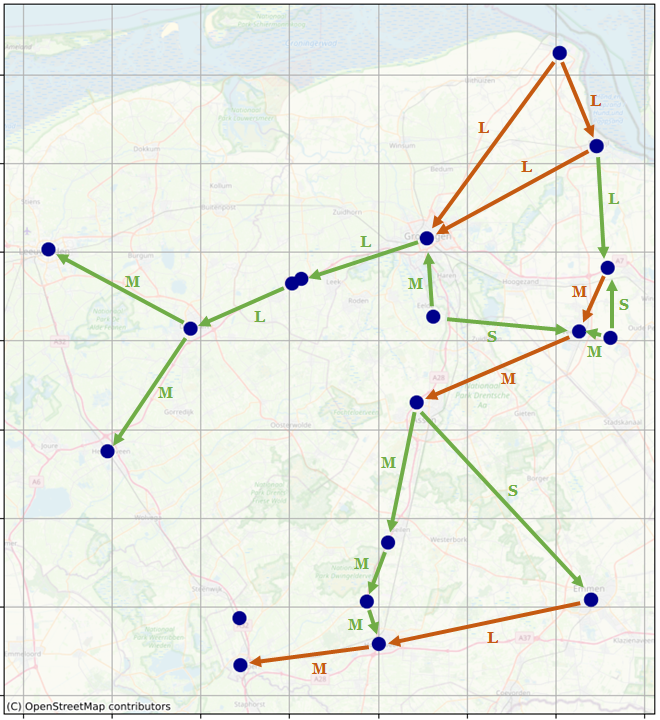}
        \caption{Decentralized production on deterministic case}
        \label{fig:case2-eev}
        \end{subfigure}
    \caption{  Comparison of final network designs (at $t=T$) for the stochastic (top row) and deterministic (bottom row) solutions across the two supply cases. Green pipelines are for hydrogen, while orange ones are for natural gas. Arc labels (S, M, L) indicate capacity class.}
    \label{fig:init}
\end{figure}

The pipeline connections for hydrogen at the latest period ($t=T$) are given in Figure \ref{fig:init} for both variants. As a comparison, we also show the deterministic (EV) solution for both variants to highlight the need to include uncertainty in decision-making for a future hydrogen economy. For clarity in representation, when at least two pipelines exist between a given pair of locations, we denote the capacity label of the largest capacity among the pipelines. Pipelines transporting green hydrogen are green and those for natural gas are orange. Note that despite the phase-out of the natural gas economy at the latest period, the graph still includes pipelines belonging to the old natural gas network. Those are the pipelines where repurposing for hydrogen use was not deemed necessary during the planning horizon. For a more detailed discussion on this behavior, see Section \ref{sect:timing}.
%that are not repurposed for hydrogen. 
%This is due to various factors such as the need for new pipelines to meet future hydrogen demands or changes in supplier locations. 

In Figure \ref{fig:case1}, we observe a pipeline network that primarily extends from Eemshaven (most upper-right node) to various destinations, mostly optimized for the shortest paths to their respective endpoints. In contrast, our analysis of the decentralized production reveals a more crowded network in Figure \ref{fig:case2}, with a higher number of pipelines near potential supply locations. The total installation capacity remains relatively consistent between both variants. The Eemshaven variant picks for pipelines with greater capacity, resulting in fewer pipelines overall. Conversely, the decentralized production variant features a larger number of pipelines, each with smaller capacities. 
In \ref{10fig}, we provide ten figures detailing the evolution of the suggested energy network over each time period for decentralized production. The majority of the network transition is observed around the center of the horizon, specifically at $t=4$.%, as this behavior is also highlighted in our analysis in Section \ref{sect:timing}. 

%In the Eemshaven variant, the starting natural gas infrastructure illustrated in Figure \ref{fig:aprime} is substantially repurposed for hydrogen transport compared to the decentralized production variant. Specifically, the conversion costs in the Eemshaven import are $70\%$ higher relative to the new pipeline construction costs than in the decentralized production. 
The Eemshaven variant repurposes more of its initial natural gas pipeline network for hydrogen as compared to the decentralized variant. Specifically, the ratio of conversion to construction costs is $70\%$ higher in Eemshaven variant than in the decentralized variant. Since the conversion costs are lower compared to the new construction costs, this leads to $17.5\%$ less total infrastructure costs on the Eemshaven variant than the decentralized variant. However, the Eemshaven variant also experiences approximately $50\%$ more penalty costs than the decentralized variant due to the greater need for truck transportation in times of peak demand. This difference arises because the extensive network of the decentralized variant, even with lower capacities, provides a more reliable setting for unexpected increases in demand. %Nevertheless, the Eemshaven variant presents economic advantages, with construction and conversion costs being $17.5\%$ lower, benefiting from the certainty in supplier location. In contrast, the decentralized production, with an emphasis on local production, is perceived as a safer alternative as it ensures the country's self-sufficiency in hydrogen production. Moreover, with its extensive pipeline coverage and more spread supply network, this variant also enhances network resilience, better preparing it to handle a broader range of unforeseen scenarios.

%To gain insights into the structure of the proposed pipeline network, w
We additionally examine the cases where we use the mean values and solve the deterministic version of the problem. For the deterministic case of the Eemshaven, the recommended policy in Figure \ref{fig:case1-eev} spends $60\%$ of the total infrastructure costs observed in its stochastic equivalent. However, due to the high penalty costs incurred in the deterministic policies, the suggested pipeline network translates to an expected overall cost increase of $20\%$. The disparity becomes more pronounced in the decentralized production variant. In the deterministic case of the decentralized production, because only $22\%$ of the total infrastructure cost of its stochastic equivalent is allocated as in Figure \ref{fig:case2-eev}, many extreme scenarios are overlooked, resulting in a $115\%$ rise in overall expected costs. This significant increase is partly attributed to certain locations being completely omitted from green hydrogen pipeline connections due to their low likelihood of being essential in the long-term perspective, reflecting the inherent uncertainties surrounding locations. As these uncertainties increase, deterministic solutions become more inadequate, often omitting pipeline connections altogether.

Our discussions of our findings with project partners have highlighted some additional points of network supply strategies. Ports like Eemshaven, for instance, are ideal locations for hydrogen production due to their abundant access to resources essential, such as land and water. The findings from our case study indicate that centralized supply already leads to lower total infrastructure costs. However, it has a higher risk of incurring penalty costs to cover peak demand. This is due to its efficient but less dense infrastructure network as compared to the decentralized production. Thus, if decision-makers proceed with the suggested option of centralized production, we recommend they concentrate on strategic plans where the effects of uncertainties can be decreased, such as by investing in sufficient storage locations.

%However, it has a risk of high costs associated with the extreme scenarios, as seen by higher penalty costs in the Eemshaven variant. Therefore, we recommend decision-makers to concentrate on decreasing uncertainties through detailed strategic planning in the energy transition, especially for centralized production case which is recommended by the experts.

In both variants, our analysis assumes sufficient supply to meet demand. Therefore, unsatisfied demand via pipelines mainly results from the consideration that it is more cost-effective to incur penalty costs than to invest in new pipeline construction, along with the associated flow costs. In situations where supply is limited, our analysis indicates that optimal solutions prioritize meeting demand at locations closer to supply points, which is evidently more cost-efficient due to reduced construction and the flow costs for closer customers. The outcome could differ, however, if the decision maker assumes varying levels of importance to customers. For instance, higher penalty costs imposed on significant customers may influence the optimal solution, potentially leading to different pipeline structures to meet their demand. %However, our analysis operates under the assumption of the customers are equally important, having a single penalty cost of unsatisfied demand per unit.

\subsection{Summary of Case Study}
In summary, 
\begin{enumerate}[leftmargin=*,labelindent=16pt]
    \item[o] We analyzed two distinct network supply strategies: one with a single supplier at Eemshaven and the other involving multiple local suppliers adding location uncertainty.
    \item[o] The Eemshaven variant displayed a simpler pattern with pipelines radiating mostly directly from Eemshaven, while the decentralized production variant featured a more crowded network with numerous smaller-capacity pipelines near potential supply locations.
    \item[o] The Eemshaven variant exhibits a higher level of integration of the initial natural gas network, driven by the optimistic assumption of location certainty and the alignment with the existing pipeline structure's flow direction. However, this led to $50\%$ more penalty costs for the need of trucks compared to the decentralized production variant due to partially neglecting peak demand scenarios.
    \item[o] While the Eemshaven variant reduces total infrastructure costs by $17.5\%$, the single supplier location offers economies of scale and easy access to resources like land and water at the port, which is important for strategic energy transition decisions. Conversely, decentralized production ensures energy independence and extensive pipeline coverage for peak demand scenarios.
    \item[o] Finally, as uncertainty increases, deterministic policies lead to insufficient pipeline infrastructures. This results in sharp increases in penalty costs, leading to an increased total expected cost of the network design. %This trend underscores the challenges of deterministic approaches under rising uncertainties.
\end{enumerate}

\section{Conclusions} \label{sect:concl}
In this paper, we address the energy transition network design problem over a finite planning horizon. The problem focuses on strategic decision-making associated with the construction and conversion of natural gas pipelines into green hydrogen pipelines, and on balancing supply and demand in a transitioning network under uncertain conditions. We provide a  novel two-stage stochastic mixed integer programming model with continuous recourse to formulate our stochastic network design problem on energy transition. In the first stage, decisions regarding the construction and conversion of pipelines are made. The second stage, including uncertainties represented by a finite number of scenarios, determines the optimal energy flow decisions and derives corresponding inventory storage policies.

Our algorithm consistently yields solutions that are near-optimal, outperforming benchmark methods in the balance of solution quality and computational time. We also use our algorithm on a case study in the Northern Netherlands, as part of the \textit{Hydrogen Energy Applications in Valley Environments for Northern Netherlands} (HEAVENN) initiative. This real-world application not only emphasizes the practical relevance of our model but also highlights its adaptability, particularly in scenarios with high uncertainties. We observe that uncertainty in future supplier locations is a key factor for network design decisions.

Several suggestions arise for future studies on the network design topic. \uh{Firstly, a valuable extension would be to explore an adaptive multi-stage recourse model \citep{kayacik2025partially} in the context of multi-period stochastic network design, which could permit corrective investment decisions at a few critical junctures. This would offer a more nuanced approach than a pure two-stage model while still respecting the real-world constraints that make a multi-stage recourse impractical.} \uh{Second, the scenario bundling strategy used in our matheuristic can be strengthened. While this study employs random assignment as a baseline, future research should explore deterministic clustering approaches—such as k-means based on demand profiles or flow patterns—to improve convergence behavior and solution quality.} Third, considering heterogeneous customers might be practical, leading to a different penalty cost per customer. In such a setting, considering the implementation of risk-adjusted performance metrics such as Conditional Value-at-Risk (CVaR) could provide a more accurate reflection of decision-makers' risk aversion, particularly for customers involving extreme losses.  \uh{Moreover, incorporating uncertainty in pipeline construction or repurposing costs could be beneficial, as costs can vary based on geographical and geological conditions. Additionally, integrating principles of land-use planning or urban development may offer insights into how the strategic positioning of hydrogen refueling stations or households affects energy costs in the network. }

Another promising future research direction involves adapting our progressive hedging-based matheuristic to other two-stage recourse models, due to its generality. Progressive hedging algorithms often exhibit late convergence when only a small subset of decision variables remains non-convergent for a substantial number of iterations \citep{watson2011progressive}. For this, our heuristic refines the progressive hedging algorithm by applying problem-specific constraints that recognize the similarities among decision variable subsets, rather than relying solely on individual variable convergence\uh{, i.e., aggregate soft fixing}. These constraints manage cases where individual decision variables may differ, yet a consistent pattern across the scenarios may be utilized to accelerate convergence while maintaining high solution quality. Similarly, the small subset of decision variables remaining non-convergent may be addressed more efficiently using state-of-the-art optimization solvers such as Gurobi, as we demonstrated with an early convergence check in our algorithm. Future studies should explore this progressive hedging framework for two-stage recourse models, by combining \uh{aggregate soft fixing, }early convergence checks, and scenario bundling, for accelerated, high-quality solutions.

%Lastly, introducing a budget constraint for each period could ensure a more evenly distributed construction cost over the periods.
\section*{Acknowledgments}
This project has received funding from the Fuel Cells and Hydrogen 2 Joint Undertaking (now Clean Hydrogen Partnership) under Grant Agreement No 875090. This Joint Undertaking receives support from the European Union's Horizon 2020 research and innovation programme, Hydrogen Europe and Hydrogen Europe Research. Albert H. Schrotenboer has received support from the Dutch Science Foundation (NWO) through grant VI.Veni.211E.043. We thank the Center for Information Technology of the University of Groningen for their support and for providing access to the Hábrók high performance computing cluster.

\newpage
\bibliographystyle{elsarticle-harv} 
\bibliography{main}

\appendix

\newpage
\section{Instance Generation} \label{instGen}

In this section, we outline the procedure for generating instances for our study. Broadly, each instance is characterized by a \textbf{graph}, associated \textbf{costs}, and instance \textbf{parameters}. We discuss these components in sequence in this section. We determine costs and parameters on the expert recommendations of the HEAVENN project, and tailored them for our computational experiments, ensuring valuable insights for readers. Parameters associated with the solution methodology are elaborated in Section \ref{sect:soln}.

For the \textbf{graph}, nodes are uniformly randomized within a $10\times10$ square area (which represents approximately $100\times100$ km square area). Subsequently, nodes are categorized based on their functional roles. A node can be designated as
\begin{enumerate}[itemsep=0pt, parsep=0pt]
    \item a supply source,
    \item a hydrogen refueling station or an analogous demand unit for natural gas,
    \item a household region requiring gas for heating,
    \item an industrial hub requiring gas as a feedstock,
    \item and an intermediate node for pipeline connections.
\end{enumerate} 

Then, nodes are assigned these roles based on the following probabilities: $\{0.2,0.3,0.2,0.1,0.2\}$, respectively, different for each commodity. However, we allow for the possibility that certain locations may assume different roles depending on the scenario, introducing an element of uncertainty particularly for the type of location for supply (1) and intermediate nodes (5). Thus, a location selected as a supply node in one scenario of an instance might serve as an intermediate node in another scenario of the same instance. We do not, however, introduce this uncertainty for demand locations (2, 3, and 4), as the demand nodes are known beforehand with high certainty, e.g. household areas. Then, we establish potential pipeline connections between nodes. Specifically, given a predefined value of $U^d$ and the distance $d_{ij}$ between nodes $i$ and $j$, if there exists at least one node $\ell$ such that $d_{ij} * U^d > d_{i\ell} + d_{\ell j}$, the candidate arcs between nodes $i$ and $j$ are all discarded. This reflects the practical observation that pipelines typically do not span long distances directly between two nodes. Instead, they traverse multiple intermediate nodes.  For our computational experiments, we set $U^d$ to $1.15$. If no such $\ell$ node exists, then we add a total of $6$ arcs to $\mathscr{A}$, three pipelines in each direction between $i$ and $j$, each with diameters of $\{30,75,120\}$.% \uh{in mm}.

Following the graph generation, we define the associated \textbf{costs}. The fixed cost for pipeline construction %\uh{(in million euros)} 
is determined using the formula from \cite{andre2013design} as $(a_0 + a_1 u_{ij} + a_2 u_{ij}^2) \times d_{ij}$ where $a$ parameters are coefficients and $u_{ij}$ represents the diameter. Based on an area-specific study \citep{kevinMSC}, we select the coefficients as $a_0= 25/6$, $a_1=5/72$, and $a_2=1/5400$, by normalizing them to our distance scheme. We implement these values for all $f_{a0}^k$. Motivated by considerations of opportunity costs and economies of scale, the fixed costs for each pipeline linearly decrease to $50\%$ at the end of the horizon. In other words, $f_{a,T-1}^k$ equals $50\%$ of $f_{a0}^k$ for every $a \in \mathscr{A}$ and $k \in \mathscr{K}$. The conversion cost, denoted as $g_{at}^k$, is set to half the fixed cost of construction, i.e., $g_{at}^k \coloneqq 0.5 f_{at}^k$ for every $a \in \mathscr{A}$, $k \in \mathscr{K}$, and $t \in \mathscr{T}$. As for the flow cost, $c_{at}^k$, we establish it as $c_{at}^k = 2.5\times d_{a}$% \uh{ per tonne}
, with $d_{a}$ being the distance of arc $a$. For the terminal period, the flow cost is adjusted to account for future flow costs in a fully hydrogen-transitioned network: $c_{a,T-1}^k$ is further multiplied by $T$. To deter penalties, we assign $h$ a value of $1$.

Lastly, we determine the problem \textbf{parameters}. The scenario probability, $p_\omega$, is uniformly set to $1/ |\Omega|$. The arc capacity, $U_a$, is adapted to our context as $U_a \coloneqq (2/3) u_a  d_a $, with $u_a$ representing the diameters from the set $\{30,75,120\}$. The start inventory values of $o_{i}^{k}$ are selected as zero for all. Node capacities, $C_{it}^{k\omega}$, are categorized into three primary standards:
\begin{enumerate}[itemsep=0pt, parsep=0pt]
    \item Region-specific high-capacity storage facilities (e.g., salt caverns for hydrogen),
    \item Hub-specific medium capacity storage facilities (e.g., at industry hubs),
    \item Short-term storage units (e.g., cylinder storage for compressed gases).
\end{enumerate}

These standards correspond to capacities of $\{10000,1000,10\}$ %\uh{in kg}
, respectively. Table \ref{table:probMat} illustrates the probability distribution of these capacities across node types for every $i \in \mathscr{N}$, $k \in \mathscr{K}$, and $\omega \in \Omega$. We assume that $C_{it}^{k\omega}$ remains the same for all $t \in \mathscr{T} \cup \{T\}$.

{
\renewcommand{\arraystretch}{0.85}
\begin{table}[h]
\centering
\caption{Probability of storage units per each location type}
\label{table:probMat}
\begin{tabular}{r|rrrrr}
Capacity      & Supplier & HRS & Household & Industry & Intermediate \\ \hline
10000 & 0.1      & 0.0 & 0.00      & 0.0      & 0.2          \\
1000  & 0.9      & 0.3 & 0.05      & 0.5      & 0.1          \\
10    & 0.0      & 0.7 & 0.95      & 0.5      & 0.7         
\end{tabular}
\end{table}
}

We further generate demand samples, $s_{it}^{k\omega}$, by considering that $\mathscr{K} = 2$ with hydrogen and natural gas represented by indices $k=0$ and $k=1$, respectively. We introduce \uh{$pr_t^\omega$} as the probability of observing $s_{it'}^{1\omega} = 0$ for $t' \in \{t+1,\dots,T\}$ for natural gas, and $s_{it'}^{0\omega} = 0$ for $t' \in \{0,\dots,t - 1\}$ for hydrogen. Essentially, \uh{$pr_t^\omega$} depicts the transition probability from natural gas to hydrogen at time $t'$. For each scenario $\omega$, we select \uh{$pr_t^\omega$} such that \uh{$pr_t^\omega$} $>>$ \uh{$pr_{t-1}^\omega$} in order to represent the expected exponentiality of the transition from natural gas to hydrogen, where $\uh{pr_0^\omega} = 0$ and $\sum_{t \in \mathscr{T}} \uh{pr_t^\omega} = 1$. Specifically, we select $\uh{pr_t^\omega}$ as;
\begin{align}
    \uh{pr_t^\omega} = \frac{\uh{r(\omega)^t}-1}{\sum_{t=0}^{T-1}  (\uh{r(\omega)^t}-1)} \label{a1}
\end{align}
where $r$ is uniformly randomized in $[1.01,1.4]$, differently per each scenario of each instance. Using the probability array \uh{$pr_t^\omega$} , we randomly assign a transition time matrix, \uh{$a_{ik}^\omega$}, indicating the period when commodity $k$ transitions on node $i$. For example, if $s_{it'}^{k\omega} = 0$ for $t' \in \{t+1,\dots,T\}$ for natural gas, then $\uh{a_{ik}^\omega} = t$, where $\mathbb{P} \{ \uh{a_{ik}^\omega} = t \} = \uh{pr_t^\omega}$. For our base scenario, the net supply levels at the beginning are represented by $s_{i0}^{k}$, structured to account for two distinct seasons: winter and summer. Winter mostly has heightened demand but diminished supply, whereas summer experiences the opposite. Each period $t$ alternates between these seasons. The initial net supply levels %\uh{in tonne} 
are given in the Table \ref{table:s0}.

{
\renewcommand{\arraystretch}{0.85}
\begin{table}[h]
\centering
\caption{Initial net supply levels of locations per season}
\label{table:s0}
\begin{tabular}{r|rrrrr}
Season & Supplier & HRS & Household & Industry & Intermediate \\ \hline
summer & 240      & -30 & -20       & -50      & 0            \\
winter & 120      & -40 & -40       & -40      & 0           
\end{tabular}
\end{table}
}

The generation of $s_{it}^{k\omega}$ is further detailed in Algorithm \ref{alg:sitkw}. In this methodology, the growth factor $\uh{G^\omega}$ represents the expected demand increase ratio on each period over the planning horizon. We introduce an additional factor, $H$, specifically for hydrogen to capture the potentially elevated growth rates due to its emergence as a new energy source. Additionally, we introduce $R$ to account for the uncertainty in short-term energy demand fluctuations. We pick values $[G_L,G_U] = [1.0,1.1]$, $[H_L,H_U] = [1.3,1.7]$, and $R=0.2$ in our computational experiments.

\begin{algorithm}[h]
    \caption{Generation of $s_{it}^{k\omega}$} \label{alg:sitkw}
    \hspace*{\algorithmicindent} \textbf{Input:} $s_{i0}^{k}$, $[G_L,G_U]$, $[H_L,H_U]$, $R$.\\
    \hspace*{\algorithmicindent} \textbf{Output:} $s_{it}^{k\omega}$
    \begin{algorithmic} [1]
        \ForEach {scenario \uh{$\omega \in \Omega$}}
        \State Randomize a probability array \uh{$pr_t^\omega$} for $t \in \mathscr{T}$.
        \State Assign transition time matrix; \uh{$a_{ik}^\omega$}, according to \uh{$pr_t^\omega$}.
        \State Randomize uniformly $\uh{G^\omega} \in[G_L,G_U]$ and $\uh{H^\omega} \in[H_L,H_U]$.
        \State For natural gas, assign $s_{it}^{k} = s_{i0}^{k} \uh{G^\omega_t}$ for each node $i$ for $t < \uh{a_{ik}^\omega}$, and $0$ for $t \geq \uh{a_{ik}^\omega}$ by the corresponding season.
        \State For hydrogen, assign $s_{it}^{k} = s_{i0}^{k} \uh{G^\omega_t H^\omega}$ for each node $i$ for $t \geq \uh{a_{ik}^\omega}$, and $0$ for $t < \uh{a_{ik}^\omega}$ by the corresponding season.
        \State For all node $i$, commodity $k$ and period $t$, select $s_{it}^{k\omega}$ in $[(1-R\frac{t}{T}) s_{it}^{k},(1+R\frac{t}{T}) s_{it}^{k\omega}]$ uniformly random.
        \EndFor 
    \end{algorithmic}
\end{algorithm}

\uh{To illustrate Algorithm \ref{alg:sitkw}, Figure \ref{fig:avg_transitions} shows the scenario-averaged net supply and demand over time for representative nodes for both commodities. The plots demonstrate the emergent gradual transition from natural gas to hydrogen% that results from our scenario generation process
. The trend shows the expected decrease in natural gas and increase in hydrogen over the planning horizon, while the oscillations, particularly visible for the natural gas, correspond to the modeled seasonality between winter and summer periods.}

\begin{figure}[H]
    \centering
    % Figure for Supply Transition
    \begin{subfigure}[b]{0.49\textwidth}
        \centering
        \begin{tikzpicture}[scale=0.85]
            % Y-Axis 1 (Left): Hydrogen Supply
            \begin{axis}[
                width=\linewidth,
                xlabel={Time Period},
                ylabel={Hydrogen Supply},
                ymin=0,
                legend pos=north west,
                y axis line style={blue},
                ylabel style={color=blue},
                yticklabel style={color=blue}
            ]
            \addplot[color=blue, mark=*, smooth] coordinates {
                (1, 0) (2, 3.03874) (3, 3.35407) (4, 7.16941) (5, 12.706) (6, 25.3186) (7, 57.4334) (8, 74.0696) (9, 93.5249) (10, 112.148) (11, 147.425) (12, 203.581) (13, 259.301) (14, 321.439) (15, 424.243)
            };
            %\addlegendentry{Hydrogen}
            \end{axis}

            % Y-Axis 2 (Right): Natural Gas Supply
            \begin{axis}[
                width=\linewidth,
                axis y line*=right,
                axis x line=none,
                ylabel={Natural Gas Supply},
                ymin=0,
                legend style={at={(0.97,0.97)},anchor=north east},
                % --- CORRECTED THIS PART ---
                color=red,
                % --------------------------
                ylabel style={color=red},
                yticklabel style={color=red},
                ylabel near ticks,
                ytick pos=right,
                yticklabel pos=right
            ]
            \addplot[color=red, mark=x, smooth] coordinates {
                (1, 240.0000) (2, 211.2293) (3, 249.2087) (4, 209.9419) (5, 241.1137) (6, 225.1514) (7, 242.2860) (8, 236.4707) (9, 240.3805) (10, 225.5964) (11, 209.1356) (12, 176.3089) (13, 133.3735) (14, 74.4721) (15, 0.0000)
            };
            %\addlegendentry{Natural Gas}
            \end{axis}
        \end{tikzpicture}
        \caption{   Examples of net supply transitions for selected supply nodes}
        \label{fig:supply_transition}
    \end{subfigure}
    \hfill
    % Figure for Demand Transition
    \begin{subfigure}[b]{0.49\textwidth}
        \centering
        \begin{tikzpicture}[scale=0.85]
            % Y-Axis 1 (Left): Hydrogen Demand
            \begin{axis}[
                width=\linewidth,
                xlabel={Time Period},
                ylabel={Hydrogen Demand},
                legend pos=south west,
                y axis line style={blue},
                ylabel style={color=blue},
                yticklabel style={color=blue}
            ]
            \addplot[color=blue, mark=*, smooth] coordinates {
                (1, 0) (2, -0.54354) (3, -0.458608) (4, -1.86758) (5, -2.05268) (6, -7.44124) (7, -9.03244) (8, -14.6088) (9, -13.988) (10, -26.6787) (11, -25.4863) (12, -51.7616) (13, -52.7028) (14, -103.798) (15, -102.357)
            };
            %\addlegendentry{Hydrogen}
            \end{axis}

            % Y-Axis 2 (Right): Natural Gas Demand
            \begin{axis}[
                width=\linewidth,
                axis y line*=right,
                axis x line=none,
                ylabel={Natural Gas Demand},
                legend style={at={(0.97,0.03)},anchor=south east},
                % --- CORRECTED THIS PART ---
                color=red,
                % --------------------------
                ylabel style={color=red},
                yticklabel style={color=red},
                ylabel near ticks,
                ytick pos=right,
                yticklabel pos=right
            ]
            \addplot[color=red, mark=x, smooth] coordinates {
                (1, -30) (2, -41.8558) (3, -31.9316) (4, -44.5277) (5, -34.0046) (6, -46.5015) (7, -35.4273) (8, -45.7513) (9, -30.5097) (10, -37.143) (11, -26.4931) (12, -29.8717) (13, -17.4303) (14, -11.3237) (15, 0)
            };
            %\addlegendentry{Natural Gas}
            \end{axis}
        \end{tikzpicture}
        \caption{   Examples of net supply transitions for selected demand nodes}
        \label{fig:demand_transition}
    \end{subfigure}
    \caption{  Illustration of the emergent gradual transition for average net supply (a) and demand (b) for $T=15$ and $W=120$. The overall trends show the shift from natural gas to hydrogen, while the oscillations reflect modeled seasonality.}
    \label{fig:avg_transitions}
\end{figure}

\newpage

\section{Benders decomposition} \label{benders}

In this appendix, we detail the Benders decomposition approach applied to our two-stage stochastic network design problem. The first stage variables are handled in the master problem, while the second stage decisions are decided in the subproblem.

\textbf{Master Problem:}
\begin{align}
    \text{min  } & \sum_{a \in \mathscr{A}} \sum_{t \in \mathscr{T}} \sum_{k \in \mathscr{K}} \left( f_{at}^k y_{at}^k  +  g_{at}^k b_{at}^k \right) + \sum_{\omega \in \Omega} p_\omega \theta_\omega \span \label{master_objective} \\
    \text{s.t  } & \sum_{k \in \mathscr{K}} \sum_{t \in \mathscr{T}} y_{at}^k \leq 1 & \forall a \in \mathscr{A},  \\
    & y_{at}^{k} \leq z_{at}^{k} & \forall a \in \mathscr{A}, k \in \mathscr{K}, t \in \mathscr{T}, \\
    & \sum_{k \in \mathscr{K}} \sum_{t' = 0}^t y_{at'}^{k} = \sum_{k \in \mathscr{K}} z_{at}^{k}  & \forall a \in \mathscr{A}, t \in \mathscr{T}, \\
    & z_{at}^{k} + \sum_{k' \in \mathscr{K} \setminus \{k\}} z_{a,t-1}^{k'} \leq 1 + b_{at}^k & \forall a \in \mathscr{A}, k \in \mathscr{K}, t \in \mathscr{T} - \{0\}, \\ 
    & \sum_{k \in \mathscr{K}} \sum_{t \in \mathscr{T}} b_{at}^k \leq 1 & \forall a \in \mathscr{A}  ,\\
    & y_{a0}^k = 1 & \forall k \in \mathscr{K}, a \in \mathscr{A}_k,  \\ 
    & y_{at}^k \in \{0,1\} & \forall a \in \mathscr{A}, k \in \mathscr{K}, t \in \mathscr{T},\\
    & z_{at}^{k}\in \{0,1\} & \forall a \in \mathscr{A}, k \in \mathscr{K}, t \in \mathscr{T},\\
    & b_{at}^{k}\in \{0,1\} & \forall a \in \mathscr{A}, k \in \mathscr{K}, t \in \mathscr{T}, \\
    & \theta_\omega \geq 0 & \forall \omega \in \Omega.
\end{align}
where $\theta_\omega$ is the corresponding decision variable for the sub problem of scenario $\omega$.

\textbf{Sub problem:}

\textit{Primal problem:}
For each $\omega \in \Omega$;
\begin{align}
    \text{min  } & \sum_{k \in \mathscr{K}} \sum_{t \in \mathscr{T}} \left(  \sum_{a \in \mathscr{A}}   c_{at}^k x_{at}^{k\omega}  + \sum_{i \in \mathscr{N}} h e_{it}^{k\omega}  \right) \span \label{sub_primal_start} \\
    \text{s.t  } & I_{i,t+1}^{k\omega} \leq I_{it}^{k\omega} + s_{it}^{k\omega} + e_{it}^{k\omega}  - \sum_{a \in \delta^-(i)} x_{at}^{k\omega} + \sum_{a \in \delta^+(i)} x_{at}^{k\omega} &  \forall i \in \mathscr{N}, k \in \mathscr{K}, t \in \mathscr{T}, \span \label{dual_var_start}\\
    & I_{it}^{k\omega} \leq C_{it}^{k\omega} \span \forall i \in \mathscr{N}, k \in \mathscr{K}, t \in \mathscr{T}\cup \{T\},  \\
    & x_{at}^{k\omega} \leq U_a  z_{at}^{k*} & \forall a \in \mathscr{A}, k \in \mathscr{K}, t \in \mathscr{T}, \\
    & I_{i0}^{k\omega} = o_{i}^{k} &  \forall i \in \mathscr{N}, k \in \mathscr{K}, \label{phiPrimalConst} \\
    & x_{at}^{k\omega} \geq 0 & \forall a \in \mathscr{A}, k \in \mathscr{K}, t \in \mathscr{T},\\
    & I_{it}^{k\omega} \geq 0 \span \forall i \in \mathscr{N}, k \in \mathscr{K}, t \in \mathscr{T}\cup \{T\}, \\
    & e_{it}^{k\omega} \geq 0 & \forall i \in \mathscr{N}, k \in \mathscr{K}, t \in \mathscr{T}, \label{sub_primal_end}
\end{align}
where $z_{at}^{k*}$ is the optimal value of the current master problem iteration.

\textit{Dual problem:}
For each $\omega \in \Omega$;
\begin{align}
    \text{max  } & \sum_{k \in \mathscr{K}}  \sum_{t \in \mathscr{T}} \left(   \sum_{i \in \mathscr{N}} \left( s_{it}^{k\omega} u_{it}^{k\omega}  +  C_{it}^{k\omega} \psi_{it}^{k\omega} \right) + \sum_{a \in \mathscr{A}}  U_a  z_{at}^{k*} v_{at}^{k\omega}  \right)  +  \sum_{k \in \mathscr{K}} \sum_{i \in \mathscr{N}} o_{i}^{k} \phi_{i}^{k\omega}\span \label{sub_dual_start}\\
    \text{s.t  } & u_{jt}^{k\omega} - u_{\ell t}^{k\omega} + v_{at}^{k\omega} \leq c_{at}^k &  \forall a \coloneqq \{j,l\} \in \mathscr{A}, k \in \mathscr{K}, t \in \mathscr{T}, \\
    & u_{i,t-1}^{k\omega} \mathbb{1} (t>0) - u_{it}^{k\omega} \mathbb{1} (t < T) + \psi_{it}^{k\omega} + \phi_{i}^{k\omega}\mathbb{1}(t=0) \leq 0&  \forall i \in \mathscr{N}, k \in \mathscr{K}, t \in \mathscr{T} \cup \{T\}, \label{dualConstPhi} \\
    & - u_{it}^{k\omega} \leq h & \forall i \in \mathscr{N}, k \in \mathscr{K}, t \in \mathscr{T}, \\
    %& u_{it}^{k\omega} \leq h_{t}^{-k} & \forall i \in \mathscr{N}, k \in \mathscr{K}, t \in \mathscr{T} \\
    & u_{it}^{k\omega} \leq 0 & \forall i \in \mathscr{N}, k \in \mathscr{K}, t \in \mathscr{T}, \\
    & \psi_{it}^{k\omega} \leq 0 & \forall i \in \mathscr{N}, k \in \mathscr{K}, t \in \mathscr{T} \cup \{T\}, \\
    & v_{at}^{k\omega} \leq 0 & \forall a \in \mathscr{A}, k \in \mathscr{K}, t \in \mathscr{T}, \\
    & \phi_{i}^{k\omega} \text{ urs} & \forall i \in \mathscr{N}, k \in \mathscr{K}, \label{phiDefinition}
\end{align}

where $u_{it}^{k\omega}$, $\psi_{it}^{k\omega}$, $v_{at}^{k\omega}$, and $\phi_{i}^{k\omega}$ are dual decision variables for Constraints \eqref{dual_var_start} - \eqref{phiPrimalConst}, respectively.

There exists no feasibility cut to be generated, since primal subproblem in Eq. \eqref{sub_primal_start} - \eqref{sub_primal_end} is always feasible independent of $z_{at}^{k*}$. However, the optimal solution of the dual model in Eq. \eqref{sub_dual_start} - \eqref{phiDefinition} defines an optimality cut for the master problem in the form of the following;
\begin{align}
    \theta_\omega \geq  \sum_{k \in \mathscr{K}}  \sum_{t \in \mathscr{T}} \left(   \sum_{i \in \mathscr{N}} \left( s_{it}^{k\omega} u_{it}^{k\omega*}  +  C_{it}^{k\omega} \psi_{it}^{k\omega*} \right) + \sum_{a \in \mathscr{A}}  U_a  z_{at}^{k} v_{at}^{k\omega*}  \right)  +  \sum_{k \in \mathscr{K}} \sum_{i \in \mathscr{N}} o_{i}^{k} \phi_{i}^{k\omega*}.
\end{align}
where $u_{it}^{k\omega*}$, $\psi_{it}^{k\omega*}$, $v_{at}^{k\omega*}$, and $\phi_{i}^{k\omega*}$ represent the optimal solution of the dual problem.

We iteratively solve the master problem and subproblems, generating one cut per scenario on each iteration. The algorithm terminates when the relative gap between the lower bound (LB) and the upper bound (UB) closes to within $1\%$, aligning with the MIP gap used in GRB. Here, LB is defined as the value of \eqref{master_objective} of the last iteration, and UB is the best solution found so far by the algorithm, as in Eq. \eqref{firstEq}. 

\uh{Initially, we explored the use of Benders decomposition as it is common for stochastic network design problems \citep{crainic2021partial,rahmaniani2018accelerating,rahmaniani2024asynchronous}. Preliminary experiments, however, indicated significant limitations in terms of computational efficiency for our problem. In order to handle with these inefficiencies, \cite{crainic2021partial} accelerate a partial Benders decomposition with strategies such as adding cuts, MIP gaps and time limits for a stochastic network design study, resulting in optimality gaps about $1-5\%$ for problems mostly sized at $N=10$ with no multi-period setting. Similar MIP gaps, and time limits are also implemented in the other Benders decomposition based network design literature \citep{rahmaniani2018accelerating, rahmaniani2024asynchronous}, resulting in sub-optimal solutions. Because of this, as the problem complexity increases, particularly with multi-period settings, it is usual to apply heuristics rather than exact methods \citep{boland2017continuous, jiang2021soft, fragkos2021decomposition}. Overall, consistent with these approaches, we believe our algorithm offers a balance of solution quality and computational time for our stochastic network design setting, while providing optimality gaps on par with those achieved with Benders decomposition. Notably, as Benders decomposition, progressive-hedging is an exact algorithm in its base version \citep{rockafellar1991scenarios}.}

\begin{comment}
{
\renewcommand{\arraystretch}{0.85}
\begin{table}[h]
\centering
\color{red}
\captionsetup{ labelfont={color=red}, textfont={color=red} }
\caption{Results of the relative performance of Benders decomposition to GRB} \label{tab:benders}
%\resizebox{\textwidth}{!}{
\begin{tabular}{rrrrrrrr} \toprule
& & & & \multicolumn{2}{l}{Gap ($\%$)} & \multicolumn{2}{l}{Time ratio} \\
\cmidrule(r){5-6}\cmidrule(r){7-8}
$N$ & $A$& $T$ & $W$ &  min   & avg  & min  & med  \\\midrule
5 & 38 & 5 & 30 & -2.49 & -0.07 & 4.97   & 25.71  \\
  &  &   & 60 & -1.08 & 0.00  & 2.17   & 17.31  \\
  &  & 8 & 30 & -0.21 & 0.23  & 6.36   & 36.35  \\
  &  &   & 60 & -1.47 & 0.05  & 1.94   & 21.88  \\
8 & 79 & 5 & 30 & -2.83 & 5.60  & 53.62 & 284.63 \\
  &    &   & 60 & -0.40 & 5.98  & 10.98 & 97.62  \\ 
  &  & 8 & 30 & 0.09  & 6.61  & 5.32   & 47.86  \\
  &  &   & 60 & -0.09 & 6.35  & 1.65   & 20.70 
\\ \bottomrule                  
\end{tabular}
%}
\end{table}
}
\end{comment}

\newpage

\section{Distribution Production Case Timeline} \label{10fig}

\begin{figure}[h]
    \centering
    \begin{subfigure}[b]{0.49\textwidth}
        \centering
        \includegraphics[scale=0.4]{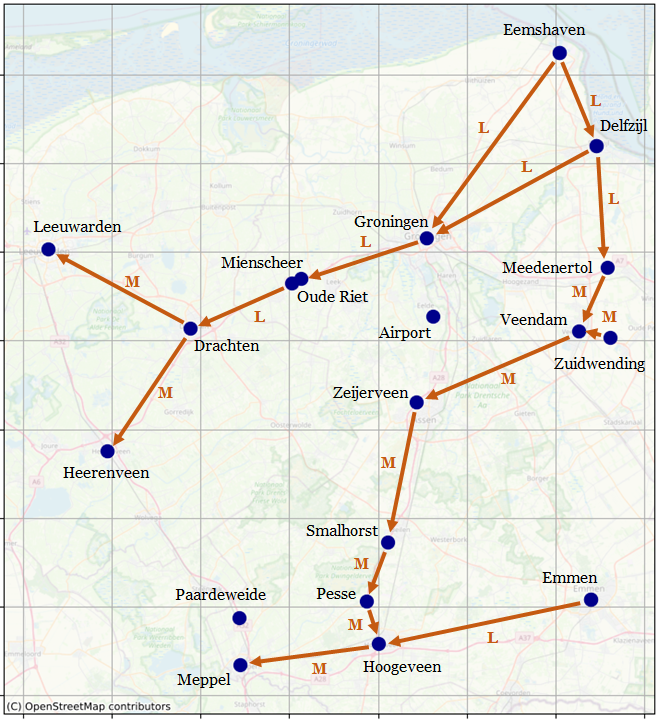}
        \caption{$t=0$}
        \label{fig:t0}
    \end{subfigure}
    \hfill
    \begin{subfigure}[b]{0.49\textwidth}
        \centering
        \includegraphics[scale=0.4]{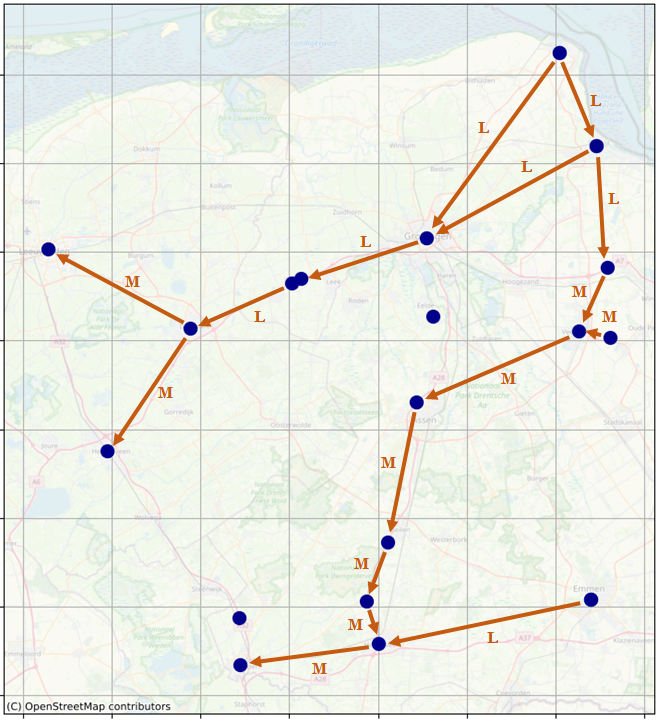}
        \caption{$t=1$}
        \label{fig:t1}
    \end{subfigure}
    \vfill
    \begin{subfigure}[b]{0.49\textwidth}
        \centering
        \includegraphics[scale=0.4]{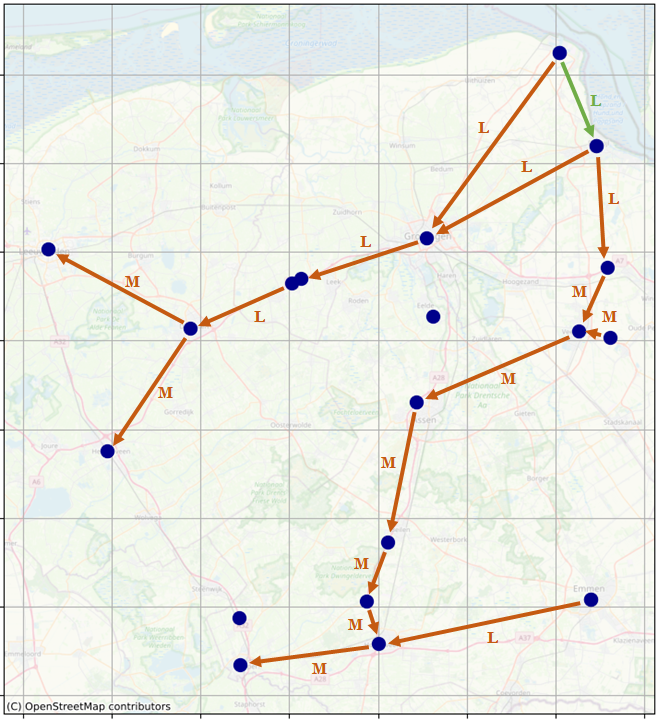}
        \caption{$t=2$}
        \label{fig:t2}
        \end{subfigure}
    \hfill
    \begin{subfigure}[b]{0.49\textwidth}
        \centering
        \includegraphics[scale=0.4]{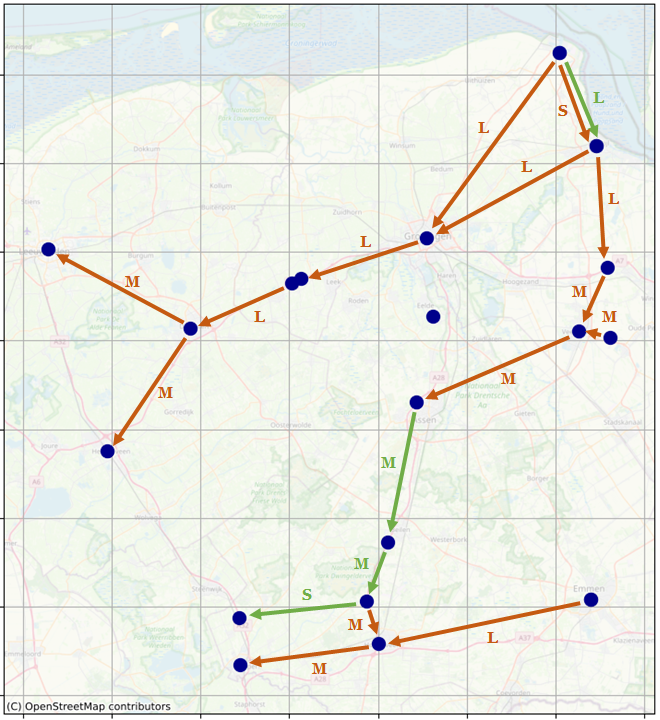}
        \caption{$t=3$}
        \label{fig:t3}
        \end{subfigure}
    \caption{The suggested pipeline network for decentralized production case at $t=0$ to $t=3$}
    \label{fig:cases1}
\end{figure}

\begin{figure}[h]
    \centering
    \begin{subfigure}[b]{0.49\textwidth}
        \centering
        \includegraphics[scale=0.40]{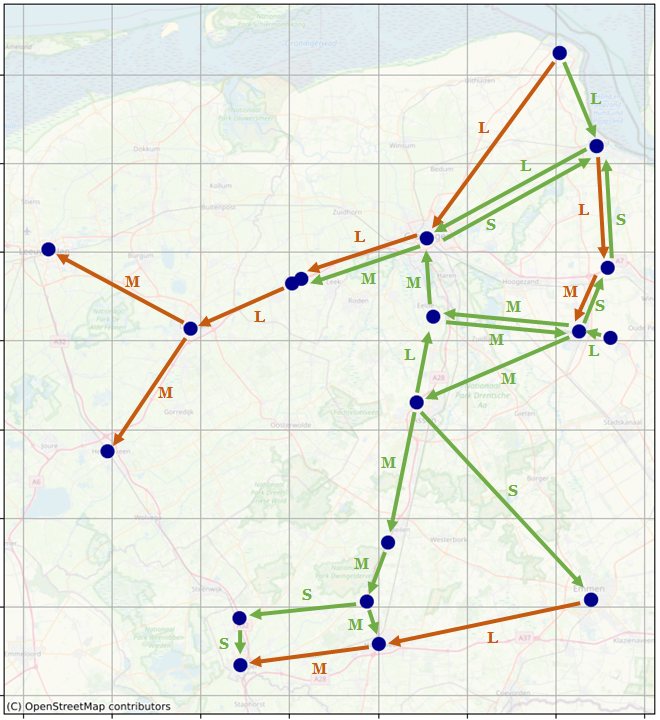}
        \caption{$t=4$}
        \label{fig:t4}
    \end{subfigure}
    \hfill
    \begin{subfigure}[b]{0.49\textwidth}
        \centering
        \includegraphics[scale=0.40]{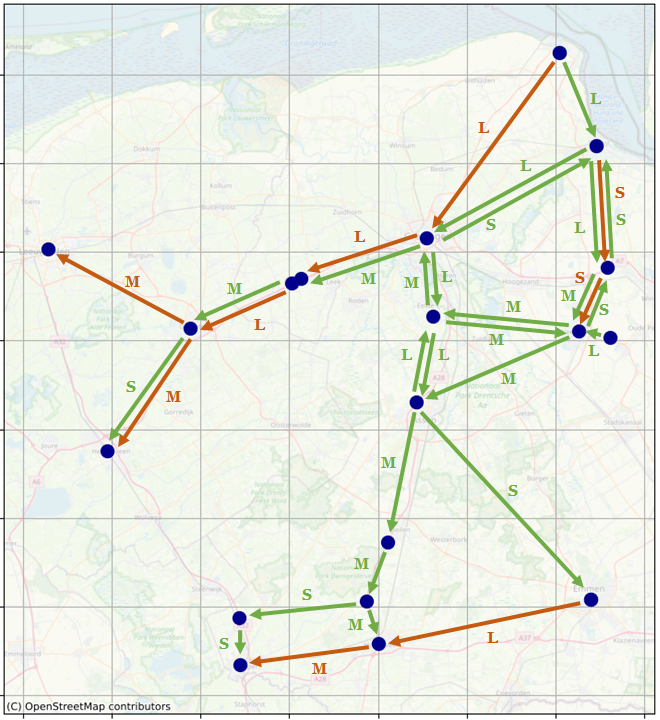}
        \caption{$t=5$}
        \label{fig:t5}
    \end{subfigure}
    \vfill
    \begin{subfigure}[b]{0.49\textwidth}
        \centering
        \includegraphics[scale=0.40]{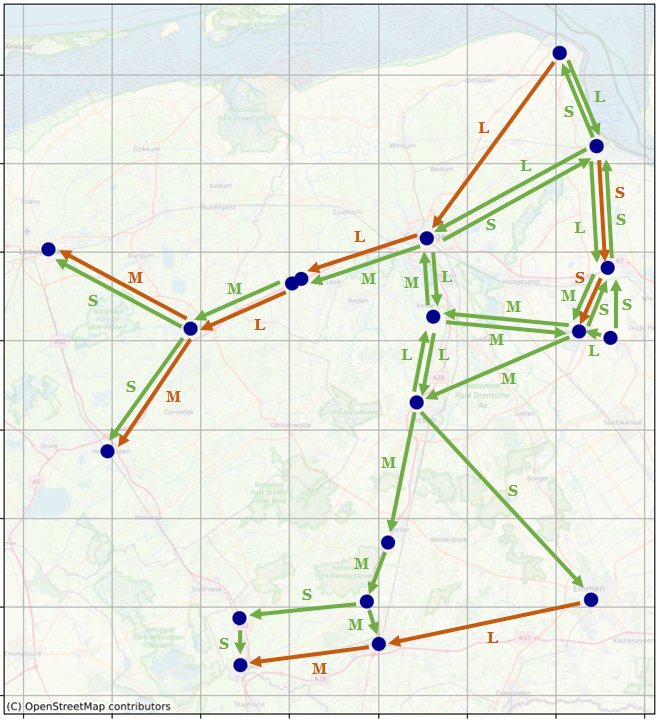}
        \caption{$t=6$}
        \label{fig:t6}
        \end{subfigure}
    \hfill
    \begin{subfigure}[b]{0.49\textwidth}
        \centering
        \includegraphics[scale=0.40]{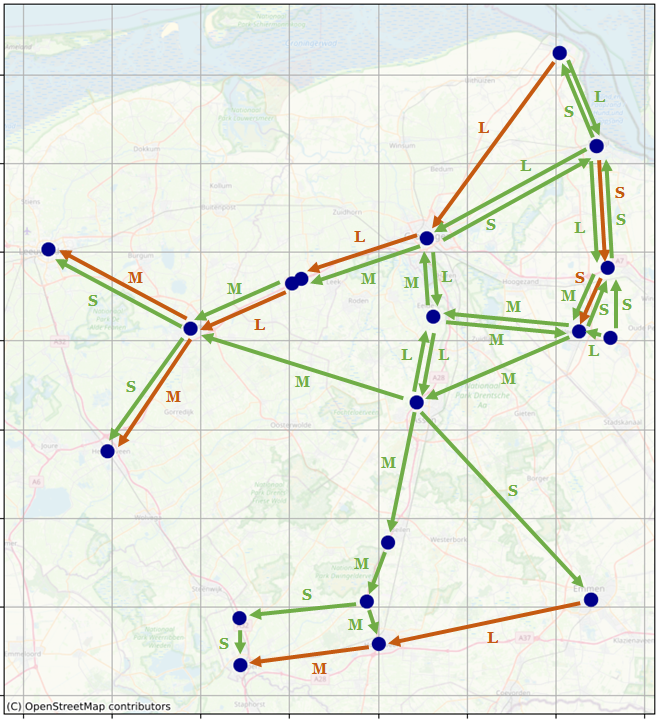}
        \caption{$t=7$}
        \label{fig:t7}
        \end{subfigure}
    \caption{The suggested pipeline network for decentralized production case at $t=4$ to $t=7$}
    \label{fig:cases2}
\end{figure}

\begin{figure}[h]
    \centering
    \begin{subfigure}[b]{0.49\textwidth}
        \centering
        \includegraphics[scale=0.40]{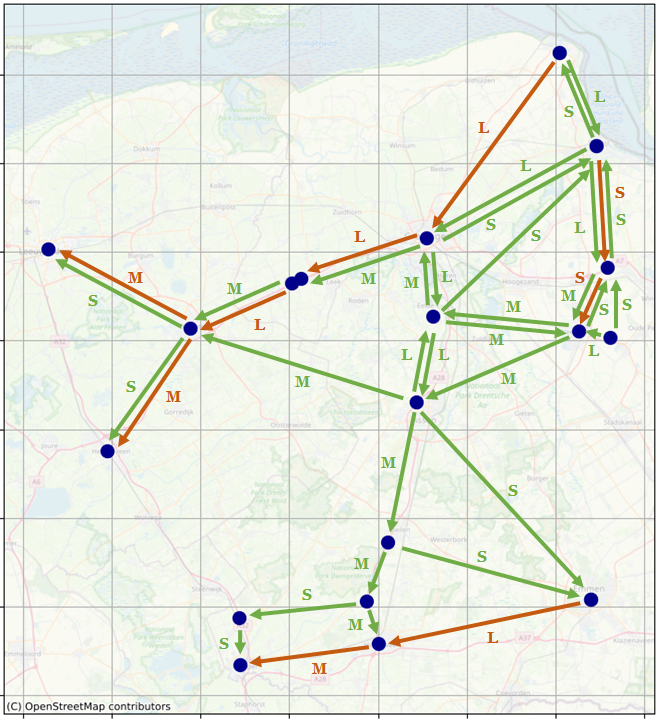}
        \caption{$t=8$}
        \label{fig:t8}
    \end{subfigure}
    \hfill
    \begin{subfigure}[b]{0.49\textwidth}
        \centering
        \includegraphics[scale=0.40]{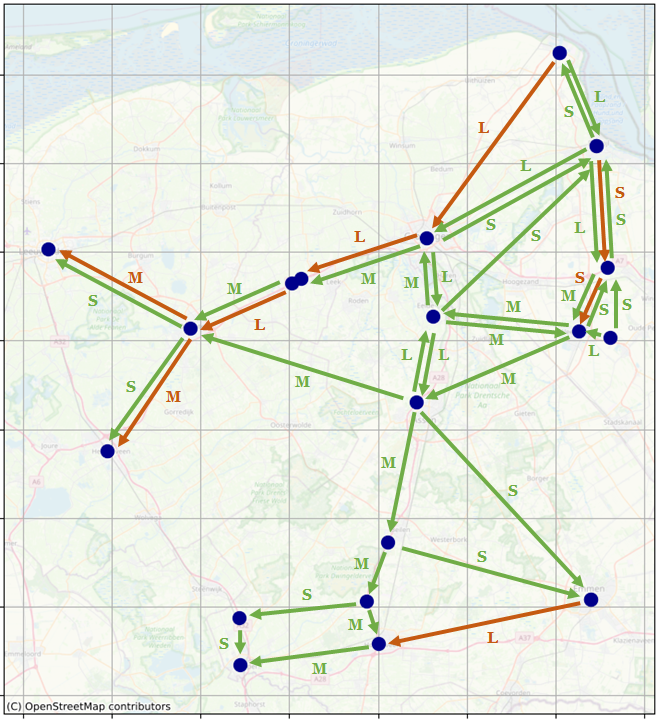}
        \caption{$t=9$}
        \label{fig:t9}
    \end{subfigure}
    \caption{The suggested pipeline network for decentralized production case at $t=8$ and $t=9$}
    \label{fig:cases3}
\end{figure}

% \bibliographystyle{elsarticle-num} 
% \bibliographystyle{elsarticle-harv} 


%\begin{thebibliography}{00}


%\bibitem{refs2}

%\end{thebibliography}
\end{document}